\documentclass[11pt]{amsart}
\usepackage[all]{xy}
\usepackage[all]{xy} \usepackage{amssymb} 
\usepackage{amsmath} 
\usepackage{graphicx} 
\usepackage[english]{babel} 
\usepackage{epsfig}
\usepackage{enumerate}
\usepackage{color}

\newenvironment{proo}{\begin{trivlist} \item{\emph{Proof.}}}
 {\hfill $\square$ \end{trivlist}}
\swapnumbers

\theoremstyle{definition}
\newtheorem{definition}[subsection]{Definition}
\newtheorem{example}[subsection]{Example}
\newtheorem{notation}[subsection]{Notation}
\newtheorem{remark}[subsection]{Remark}
\theoremstyle{plain}
\newtheorem{theorem}[subsection]{Theorem}
\newtheorem{proposition}[subsection]{Proposition}
\newtheorem{lemma}[subsection]{Lemma}
\newtheorem{corollary}[subsection]{Corollary}

\def\t{\times}
\def\ot{\otimes}
\def\K{\mathbb K}

\def\B{\mathcal B}
\def\C{\mathcal C}
\def\Irr{\mbox {\bf Irr}_n^r}

\def\ST{{\mbox {\bf ST}_n^r}}
\def\st{\mathcal {ST}\!}

\def\Br{{\mbox {\bf Br}}}
\def\Prod{{\mbox {\bf Prod}}}

\def\NN{{\mathbb{N}}}

 \newcommand{\commentc}[1]{}

\begin{document}

\author[E. Burgunder, P.-L. Curien, M. Ronco]{Emily Burgunder, Pierre-Louis Curien, Mar\'ia Ronco}
\address{EB: , Universit\'e Paul Sabatier\\
Institut de Math\'ematiques de Toulouse\\
118 route de Narbonne\\
F-31062 Toulouse Cedex 9 France}
\email{burgunder@math.univ-toulouse.fr}
\address{PLC: Laboratoire PPS, UMR 7126 CNRS\\
F-31062 Paris Cedex 9 , France}
\email{curien@pps.univ-paris-diderot.fr}
\address{MOR: IMAFI, Universidad de Talca\\ Campus Norte, Avda. Lircay s/n\\ Talca, Chile}
\email{maria.ronco@inst-mat.utalca.cl}

\title{Free algebraic structures on the permutohedra}
\thanks{Our joint work was partially supported by the Projects MathAmSud 13Math-05-LAIS and the ANR CATHRE ANR-13-BS02-0005-02. The third author's work is partially supported by FONDECYT Project 1130939. The ANR Project ANR-11-BS01-002 HOGT partially funded the joint work of the first and third authors during January, 2015.}


\begin{abstract} 
Tridendriform algebras are a type of associative algebras, introduced independently by F. Chapoton and by J.-L. Loday and the third author, in order to describe operads related to the Stasheff polytopes.  The vector space $\st$ spanned by the faces of permutohedra has a natural structure of tridendriform bialgebra, we prove that it is free as a tridendriform algebra and exhibit a basis. Our result implies that the subspace of primitive elements of the coalgebra $\st$, equipped with the coboundary map of permutohedra, is a free cacti algebra.
\end{abstract}

\maketitle

\section*{Introduction} \label{section:introduction}  

The graded vector space spanned by the set of planar rooted trees has a rich algebraic structure, coming from the different set-theoretical operations which can be performed on trees. Furthermore, when we consider the set of planar rooted trees ${\mathcal T}_n$ with a fixed number $n$ of leaves, it has a natural structure of partially ordered set whose geometric realization is a polytope of dimension $n-1$, the Stasheff polytope.

Tridendriform algebras were defined independently by F. Chapoton in \cite{Chapoton} and J.-L. Loday and the third author in \cite{Loday-Ronco}, in order to generalize the notion of dendriform algebra introduced by J.-L. Loday in \cite{Loday}, and to get a non-symmetric operad structure described in terms of  the faces of the Stasheff polytopes. 
The definitions are similar, even if they do not coincide, Chapoton\rq s one is the graded version of Loday-Ronco\rq s tridendriform version. In previous work of the first and third authors  \cite{Burgunder-Ronco},  the two notions were described in the same framework, by adding a parameter $q$: Chapoton\rq s operad coincides with the notion of $0$-tridendriform, while one recovers the original definition of tridendriform algebra of Loday and Ronco for $q=1$.

Tridendriform algebras are a particular type of non-unital associative algebras, where the associative product is the sum of certain binary operations. Many examples of associative algebras arising from combinatorial Hopf algebras, as the bialgebra of surjective maps (see \cite{Novelli-Thibon} and \cite{Palacios-Ronco}), the bialgebra of parking functions (see \cite{Novelli-Thibon}) and the bialgebra of multipermutations (see \cite{Lam-Pylyavskyy}), come from $1$-tridendriform structures. However, the $0$-tridendriform version has the advantage of giving the right definition for working in the graded differential case, as shown by F. Chapoton in \cite{Chapoton}, who gave a version of differential graded tridendriform operad which is described by the Stasheff polytope.

In  \cite{Burgunder-Ronco}, the first and the third authors also defined  $q$-Gerstenhaber-Voronov algebras as associative algebras equipped with an additional brace structure,  satisfying certain relations. When $q=0$, these algebras may be considered a non-differential version of the operad defined in \cite{Gerstenhaber-Voronov}. They  proved that the category of conilpotent $q$-tridendrifrom bialgebras is equivalent to the category of  $q$-Gerstenhaber-Voronov algebras, via the functor which associates to any coalgebra the subspace of its primitive elements.

As a direct consequence of this result, we get that  proving the freeness of a tridendriform algebra $A$ is equivalent to proving that the subspace of its primitive elements ${\mbox {Prim}(A)}$ is free as a Gerstenhaber-Voronov algebra, when $A$ is a tridendriform bialgebra. 
\medskip

The work presented here deals with the dendriform and tridendriform bialgebra structures defined on the space spanned by the faces of permutohedra. The permutohedron of dimension $n-1$ is a regular polytope which is the geometric realization of the Coxeter poset of the symmetric group $S_n$. Its faces of dimension $r$ are described by the surjective maps from $\{ 1,\dots ,n\}$ to $\{1,\dots , n-r\}$, for $0\leq r\leq n-1$.  The vector space spanned by the faces of all permutohedra has natural structures of dendriform and $q$-tridendriform bialgebras, denoted $\st_D$ and $\st_{qT}$ respectively. 
\medskip

The goal of our work is twofold: we construct a basis of $\st_D$ as a free dendriform algebra, and  a basis of $\st_{qT}$ as a free $q$-tridendriform algebra. 
We show  that the subspace ${\mbox {Prim}(\st)}$ of primitive elements of the coalgebra $\st$ is a free brace algebra for the brace structure induced by $\st _D$. Then we proceed and construct a basis ${\mathcal B}$
of ${\mbox {Prim}(\st)}$  as a free   $q$-Gerstenhaber-Voronov algebra. The last result implies that:\begin{enumerate}
\item ${\mathcal B}$ is a basis of the free $q$-tridendriform algebra $\st _{qT}$, for all $q$;
\item looking at the graded differential case, that ${\mbox {Prim}(\st)}$ is the free cacti algebra spanned by ${\mathcal B}$.\end{enumerate}

 In a recent work V. Vong, see \cite{Vong} describes combinatorial methods to study the freeness of some algebras over regular operads. Our method is essentially different, our proof relies on the following  outline: \begin{enumerate}
 \item the vector space ${\mbox {Prim}(\st)}$ is isomorphic to the vector space $\K[{\mbox {\bf Irr}}]$, spanned by the irreducible elements of ${\mbox {\bf ST}}$ for the concatenation product $\t$
 \item there exist surprisingly simple ways to define a free brace structure, respectively a free $q$-Gerstenhaber-Voronov algebra structure, on the space $\K[{\mbox {\bf Irr}}]$,
 \item the brace algebra ${\mbox {Prim}(\st_D)}$ is isomorphic to $\K[{\mbox {\bf Irr}}]$ with its brace algebra structure, while $\K[{\mbox {\bf Irr}}]$ with its $q$-Gerstenhaber-Voronov algebra structure is isomorphic to ${\mbox {Prim}(\st_{qT})}$.\end{enumerate}
 
 We hope that this type of process will provide a standard method to define new brace, $q$-Gerstenhaber-Voronov and cacti structures related to combinatorial Hopf algebras.
 \medskip
 
The paper is composed as follows: the first two sections recall basic results on coalgebras, as well as the definitions of (tri)dendriform algebras and Gerstenhaber-Voronov algebras, and the main results about tridendriform bialgebras which we use hereinafter.

In section $3$ we give the basic definitions and constructions on surjective maps needed in later sections.

Section $4$ describes the coalgebra structure of $\st$, as well as a projection from $\st$ onto the space of its primitive elements.

In section $5$ we describe the dendriform structure of $\st_D$, we define a free brace algebra structure on the space spanned by the set of irreducible surjections and we prove that it is isomorphic to the brace algebra ${\mbox {Prim}(\st_D)}$. Finally, in sections $6$ through $8$, we prove a similar result for the $q$-tridendriform algebra $\st_{qT}$.

\subsection*{Acknowledgments}  We are grateful to Jean-Yves Thibon and Jean-\linebreak Christophe Novelli for their permanent readiness to answer our questions about combinatorial aspects of Hopf algebras, and to Muriel Livernet for stimulating discussions on the cacti operad. We want also to express our gratitude to the Issac Newton Institut for Mathematics, where our joint work began during  the Programme Grothendieck-Teichm\"uller Groups, Deformations and Operads (2013).


\subsection*{Notation}  All the vector spaces considered in the present work are over $\K$, where $\K$ is a field. For any set $X$, we denote by $\K[X]$ the vector space spanned by $X$. For any $\K$-vector space $V$, we denote by $V^+:=\K\oplus V$ the augmented vector space. 


\section{Coalgebras}
\medskip

We recall the definition of coalgebra, and introduce the notation and basic results that  we need in the rest of the work.

\begin{definition} \label{definitioncoalgebra} {\rm A {\it coalgebra} over $\K$ is a vector space $C$ equipped with a linear map $\Delta : C\longrightarrow C\ot C$ which satisfies 
the coassociativity condition:
$$(\Delta \otimes Id_C)\circ \Delta = (Id_C\otimes \Delta)\circ \Delta .$$}\end{definition}

An augmentation of a coalgebra $(C, \Delta )$ is a linear map $\epsilon : C\longrightarrow \K$ such that $\cdot \circ (\epsilon \ot Id_{C}) = Id_C = \cdot \circ (Id_C\ot \epsilon )$, where $\cdot $ denotes the action of $\K$ on $C$. 

A unit of a coalgebra $(C, \Delta , \epsilon )$ is a coalgebra map $\iota : \K \longrightarrow C$ such that $\epsilon \circ \iota = Id_{\K}$, where the coalgebra structure of $\K$ is given by $\Delta _{\K}(1_{\K}) = 1_{\K}\ot 1_{\K}$.
\medskip

\begin{definition} \label{definitionprimitive} {\rm Let  $(C, \Delta)$ be a unital augmented coalgebra. An element $c\in C$ is {\it primitive} if $\Delta (c) = c\otimes 1_{\K} + 1_{\K}\otimes c$. The {\it reduced coproduct} on $C$ is defined as 
$${\overline {\Delta }}(c) := \Delta (c) - 1_{\K}\ot c - c\ot 1_{\K}.$$}\end{definition}

We denote by ${\mbox {Prim}(C)}$ the subspace of primitive elements of $C$.

It is immediate to verify that the coassociativity of $\Delta $ implies that ${\overline {\Delta}}$ is coassociative, too. We define ${\overline {\Delta }}^i: C\longrightarrow C^{\ot i}$ recursively by:\begin{enumerate}
\item ${\overline {\Delta }}^1 := Id_C$ is the identity of $C$,
\item ${\overline {\Delta }}^i:= (Id_{C^{\ot i-1}}\ot {\overline {\Delta }})\circ {\overline {\Delta }}^{i-1}$, for $i\geq 2$.\end{enumerate}

\begin{definition} \label{definitionconilpotent} A coassociative counital and unital coalgebra $(C, \Delta)$ is called {\it conilpotent} if for all $c\in {\overline C}$, there exists $n\in \NN$ such that ${\overline {\Delta }}^m(c) = 0$, for all $m\geq n$.\end{definition}
\medskip

\begin{example} \label{cotensor} {\rm Let $V$ be a $\K$-vector space. The vector space $T(V) := \bigoplus _{n\geq 1}V^{\ot n}$, where $V^{\ot n}$ denotes the tensor product $V\ot V\ot \dots \ot V$ of $V$ $n$-times, equipped with the deconcatenation coproduct:
$$\Delta ^c (v_1\ot \dots \ot v_n) = \sum _{i=1}^{n-1} (v_1\ot \dots \ot v_i)\ot (v_{i+1}\ot \dots \ot v_n),$$
is a coalgebra. 
We denote it $T^c(V)$, and call it the {\it cotensor coalgebra} over $V$. }\end{example}

Note that $T^c(V)^+$ is a unital augmented conilpotent coalgebra.

\bigskip

\bigskip

\section{Dendriform and tridendriform bialgebras}
\medskip

We recall basic results on dendriform and tridendriform bialgebras (see  \cite{Loday}, \cite{Chapoton} and \cite{Loday-Ronco}), and Gerstenhaber-Voronov algebras (see \cite{Gerstenhaber-Voronov}). 
We also describe the main results of \cite{Burgunder-Ronco}.

\begin{definition}\label{def:TriDend} {\rm Let $A$ be a vector space over $\K$.\begin{enumerate}
\item A {\it dendriform algebra} (see \cite{Loday}) structure on $A$ is a pair of binary products $\prec:A\otimes A\to A$ and $\succ:A\otimes A\to A$,
satisfying that:
\begin{enumerate}
\item $(a\prec b)\prec c = a\prec(b\prec c + b\succ c)$,
\item $(a\succ b)\prec c = a\succ(b\prec c)$,
\item $(a\prec b+ a\succ b)\succ c = a\succ(b\succ c)$.
\end{enumerate} 
\item For any $q\in \K$, a {\it $q$-tridendriform algebra} structure on $A$ is given by three binary operations $\prec:A\otimes A\to A$, $\cdot:A\otimes A\to A$ and $\succ:A\otimes A\to A$, which satisfy the following relations:
\begin{enumerate}
\item $(a\prec b)\prec c=a\prec(b\prec c + b\succ c+ q\ b \cdot c)$,
\item $(a\succ b)\prec c=a\succ(b\prec c)$,
\item $(a\prec b+ a\succ b+  q\  a\cdot b)\succ c=a\succ(b\succ c)$,
\item $(a\cdot b)\cdot c=a\cdot(b\cdot c)$,
\item $(a\succ b)\cdot c=a\succ(b\cdot c)$,
\item $(a\prec b)\cdot c=a\cdot (b\succ c)$,
\item $(a\cdot b)\prec c=a\cdot (b\prec c)$.
\end{enumerate} \end{enumerate}}
\end{definition}
\medskip

If $(A, \prec ,\cdot, \succ )$ is a $q$-tridendriform algebra, then the space $A$ equipped with the binary operations $\prec $ and $\succcurlyeq := \ q\ \cdot +\ \succ  $ is a dendriform algebra. On the other hand, for any dendriform algebra, the operation $*  = \ \succ +\ \prec $ is associative. So, dendriform and $q$-tridendriform algebras are particular cases of non-unital associative algebras.
\bigskip

\begin{definition} \label{bidendtri} A {\it dendriform bialgebra} over $\K$ is a dendriform algebra $(H, \succ ,\prec )$ equipped with a coassociative coproduct  $\Delta :H^+\longrightarrow H^+\otimes H^+$ and a counit $\epsilon : H^+\longrightarrow \K$ satisfying the following conditions:
\begin{enumerate}
\item $(\epsilon \otimes Id)\circ \Delta (x)=1_{\K}\otimes x$ and $(Id\otimes \epsilon)\circ \Delta (x)=x\otimes 1_{\K}$,
\item $\Delta (x\succ  y):=\sum (x_{(1)}*y_{(1)})\otimes (x_{(2)}\succ  y_{(2)})$, 
\item $\Delta (x\prec  y):=\sum (x_{(1)}*y_{(1)})\otimes (x_{(2)}\prec  y_{(2)})$, \end{enumerate}
for all $x,y\in H$ where $*\  =\ \succ \ +\ \prec$, $\Delta (x)=\sum x_{(1)}\otimes x_{(2)}$, and by convention:\begin{itemize}
 \item $(x*y)\otimes (1_{\K}\succ 1_{\K}):=(x\succ  y)\otimes 1_{\K}$,  
 \item $(x*y)\otimes (1_{\K}\prec 1_{\K}):=(x\prec  y)\otimes 1_{\K}$, for $x,y\in H$.\end{itemize}
A \emph{$q$-tridendriform bialgebra} is a $q$-tridendriform algebra $H$ with a coproduct $\Delta $ such that:
 \begin{enumerate}
\item $(H,\succcurlyeq , \prec )$ is a dendriform bialgebra,
\item $\Delta (x\cdot  y):=\sum (x_{(1)}*y_{(1)})\otimes (x_{(2)}\cdot  y_{(2)})$,
\end{enumerate}
where $(x*y)\otimes (1_{\K}\cdot 1_{\K}):=(x\cdot  y)\otimes 1_{\K}$.
\end{definition}
\medskip

We observe that if $(H, \succ , \cdot , \prec , \Delta )$ is a $q$-tridendriform bialgebra, then  $(H, \succcurlyeq, \prec , \Delta )$ is a dendriform bialgebra and $(H^+, *, \Delta ^+)$ is a bialgebra in the usual sense.
\bigskip

For any bialgebra $H$, the subspace ${\mbox {Prim}(H)}$ has a natural structure of Lie algebra, but in the case of dendriform and $q$-tridendriform bialgebras, the Lie bracket comes from finer structures. 
\medskip

\begin{definition}\label{def:Gerstenhaber-Voronov} (see \cite{Gerstenhaber-Voronov})
{\rm A {\it brace} algebra is a vector space $B$ equipped with $n+1$-ary operations $M_{1n}:B\ot B^{\ot n}\longrightarrow B$, for $n\geq 0$, which satisfy the following conditions:
\begin{enumerate} \item $M_{10}= Id_B$,
\item $M_{1m}(M_{1n}(x;y_1,\dots ,y_n);z_1,\dots ,z_m)=$
$$\sum_{0\leq i_1\leq j_1\leq \dots \leq j_n\leq m} M_{1r}(x;z_1,\dots ,z_{i_1},M_{1l_1}(y_1;\dots ,z_{j_1}),\dots ,M_{1l_n}(y_n;\dots , z_{j_n}),\dots ,z_m),$$
for $x,y_1,\dots ,y_n,z_1,\dots ,z_m\in B$, where $l_k=j_k-i_k$, for $1\leq k\leq n$, and $r=\sum _{k=1}^ni_k+m-j_n+n$.\end{enumerate}
\medskip 

A $q$-\emph{Gerstenhaber-Voronov} algebra, ${\mbox {GV}_q}$ algebra for short, is 
a vector space $B$ endowed with a brace structure given by operations $M_{1n}$ and an associative product $\cdot$, 
satisfying the distributive relation:
$$\displaylines {
M_{1n}(x\cdot y;z_1,\dots ,z_n)=\hfill\cr
\hfill\sum _{0\leq i\leq j\leq n} q^{j-i}M_{1i}(x;z_1,\dots ,z_i)\cdot
z_{i+1}\cdot \ldots \cdot z_j\cdot M_{1(n-j)}(y;z_{j+1},\dots ,z_n),\cr }$$
for $x,y,z_1,\dots ,z_n\in B$, where for $q =0$ we fix that $q^j = 0$ if $j\geq 1$ and $q^0 = 1$.
}
\end{definition}
\medskip

Even if brace algebras and ${\mbox {GV}}$ algebras have an infinite number of operations and seem much more complicated than dendriform and tridendriform algebras, the type of relations that these operations satisfy allow us to give an easy recursive formula for linear bases of the free objects of both theories. 
\medskip

Given a set $X$, let ${\mathbb M}(X)= \bigcup _{n\geq 1}{\mathbb M}_n(X)$ be subset of the free brace algebra ${\mbox {Br}(X)}$ over $X$, defined recursively by:\begin{enumerate}
\item $M_0(X): = X$, 
\item $M_1(X): = \{ {\mathbb M}_{1m}(x; y_1,\dots ,y_m)\mid x, y_1,\dots ,y_m\in X, m\geq 1\}$,
\item $M_n(X): = \{ M_{1m}(x; y_1,\dots ,y_m)\mid x\in X\ {\rm and}\ y_i\in M_{j_i}(X)\ {\rm for}\ 1\leq i\leq m,\ {\rm such\ that}\  {\displaystyle \sum _{i=1}^mj_i = n-1}\}$.\end{enumerate}
\medskip

For instance, for $x, y_1,y_2, y_3 ,y_4\in X$, 
\begin{align} 
M_{14}(x; y_1,\dots ,y_4)& \in {\mathbb M}_1(X),\\
M_{13}(x; y_1, M_{11}(y_2;y_3), y_4) & \in  {\mathbb M}_2(X),\\
M_{12}(x; M_{11}(y_1;y_2),M_{11}(y_3;y_4))& \in  {\mathbb M}_3(X).
\end{align}
\bigskip

The following Lemma is an immediate consequence of Definition \ref{def:Gerstenhaber-Voronov}.

\begin{lemma} \label{bracebasis} For a set $X$, the set $M(X):=\bigcup _{n\geq 1}M_n(X)$ is a basis, as a vector space, of the free brace algebra ${\mbox {Br}(X)}$ spanned by $X$.\end{lemma}
\bigskip

In a similar way, for any set $X$, we define the subset $G(X)= \bigcup _{n\geq 1} G_n(X)$ of the free ${\mbox {GV}_q}$ algebra ${\mbox {GV}_q(X)}$ over $X$, recursively as:\begin{enumerate}
\item $G_0(X): = X$, 
\item $G_1(X): = $
$$\{ y= M_{1m}(x; y_1,\dots ,y_m)\ {\rm or}\ y = y_1\cdot \ldots \cdot y_m\mid x, y_1,\dots ,y_m\in X, m\geq 1\},$$
\item For $n\geq 2$, the set $G_n(X)$ is the disjoint union of the subsets 
$$\{ y =M_{1m}(x; y_1,\dots ,y_m)\  \mid x\in X\ {\rm and}\ y_i\in G_{j_i}(X)\ {\rm for}\ 1\leq i\leq m,\  \sum _{i=1}^mj_i = n-1\},$$ and $$\{ y = y_1\cdot \ldots \cdot y_m\ \mid \ y_i\in G_{k_i}(X)\ {\rm for}\ 1\leq i\leq m,\  \sum _{i=1}^mk_i = n\}.$$\end{enumerate}
\medskip

As in the case of free brace algebras,  Definition \ref{def:Gerstenhaber-Voronov} implies the following result.

\begin{lemma} \label{GVbasis} For a set $X$, the set $G(X):=\bigcup _{n\geq 1}G_n(X)$ is a basis, as a vector space, of the free ${\mbox {GV}_q}$ algebra spanned by $X$, denoted ${\mbox {GV}_q(X)}$ .\end{lemma}
\bigskip

\bigskip

\begin{notation} \label{notationomega} Let $(A, \prec , \succ )$ be a dendriform algebra. For a family of elements $y_1,\dots ,y_r$ in $A$, let $\omega^{\prec}(y_1,\dots ,y_r)$ and $\omega ^{\succ}(y_1,\dots ,y_r)$ be the following elements of $A$:\begin{enumerate}
\item $\omega ^{\prec}(y_1,\dots ,y_r):=y_1\prec (y_2\prec (\dots \prec (y_{r-1}\prec y_r)))$,
\item $\omega ^{\succ}(y_1,\dots ,y_r):=(((y_1\succ y_2)\succ y_3) \succ \dots )\succ y_r$.\end{enumerate}\end{notation}
\medskip

There exists a functor from the category of dendriform algebras to the category of brace algebras. 

\begin{definition} \label{Bracefromdend} Let $(A, \succ , \prec )$ be a dendriform algebra. Define operations $M_{1n}: A^{\ot (n+1)}\longrightarrow A$ as follows: 
$$M_{1n}(x;y_1,\dots ,y_n):=\sum _{r=0}^n(-1)^{n-i}\omega ^{\prec}(y_1,\dots ,y_r)\succ x\prec \omega ^{\succ}(y_{r+1},\dots ,y_n),$$
for $n\geq 1$.\end{definition}
\medskip

In \cite{Ronco}, we proved that for any dendriform algebra $(A, \succ , \prec )$, the underlying vector space $A$ with the $n+1$-ary operations $M_{1n}$ is a brace algebra. In the same  work we showed the following result:

\begin{proposition} \label{propositionbrace} For any dendriform bialgebra $H$ the subspace ${\mbox {Prim}}(H)$ is closed under the brace operations.The linear map $\varphi: T^c({\mbox {Prim}}(H))\longrightarrow H$, given by:
$$\varphi (y_1\ot \dots \ot y_r) := \omega ^{\succ}(y_1,\dots ,y_r),$$
for $y_1,\dots ,y_r \in {\mbox {Prim}}(H)$ and $r\geq 1$, is a coalgebra epimorphism.\end{proposition}
\medskip

The following Theorem is proved in \cite{Ronco}.

\begin{theorem} \label{dendstruct} For any set $X$, the free dendriform algebra ${\mbox {Dend}(X)}$ over $X$ has a natural structure of bialgebra. There exists a functor ${\mathcal U}_{dend}$ from the category ${\mbox {Brace}_{\K}}$ of brace algebras to the category of ${\mbox {BiDend}_{\K}}$ of dendriform bialgebras, left adjoint to ${\mbox {Prim}}$, satisfying that any conilpotent dendriform bialgebra $H$ is isomorphic to ${\mathcal U}_{dend}({\mbox {Prim}(H)})$.\end{theorem}
\bigskip

In \cite{Burgunder-Ronco}, we proved that the functor from the category ${\mbox {Dend}_{\K}}$, of dendriform algebras over $\K$, into the category ${\mbox {Brace}_{\K}}$, which maps $(A, \succ , \prec )$ into $(A, \{M_{1n}\}_{n\geq 1})$, composed with the functor from ${\mbox {Tridend}_{q\K}}$ to ${\mbox {Dend}_{\K}}$ factorizes through the category of ${\mbox {GV}_q}$ algebras.  That is $(A, \{M_{1n}\}_{n\geq 1}, \cdot )$ is a ${\mbox {GV}_q}$ algebra, for all $q$-tridendriform algebra $(A, \succ ,\cdot , \prec)$ and we get 
$$\begin{array}[c]{ccc}
 {\mbox {Tridend}_{q\K}}& {\rightarrow}&{\mbox {Dend}_{\K}}\\
\downarrow&&\downarrow\\
 {\mbox {GV}_{q\K}}&{\rightarrow}& {\mbox {Brace}_{\K}}\end{array}$$

Moreover, if $(H, \succ ,\cdot , \prec, \Delta )$ is a $q$-tridendriform bialgebra, then ${\mbox {Prim}}(H)$ is closed under the brace operations $M_{1n}$ and the associative product $\cdot $. Finally, we got the tridendriform version of Theorem \ref{dendstruct}:

\begin{theorem}\label{idempotents} Let $X$ be a set, the free $q$-tridendriform algebra ${\mbox {Tridend}_q(X)}$ over $X$ is isomorphim, as a coalgebra, to the cotensor coalgebra $T^c({\mbox {GV}_q(X)})$, where ${\mbox {GV}_q(X)}$ denotes the free ${\mbox {GV}_q}$ algebra spanned by $X$. The functor 
${\mbox {Prim}}: {\mbox {coBiTridend}_{q\K}}\longrightarrow {\mbox {GV}_{q\K}}$ is an equivalence of categories, where ${\mbox {coBiTridend}_{q\K}}$ denotes the category of conilpotent tridendriform bialgebras.\end{theorem}
\bigskip

\bigskip

\section{Permutations and surjective maps}
\medskip

We develop first some basic definitions and notations about surjective maps and shuffles.

For any positive integer $n\in \NN$, let $[n]$ be the finite set $\{1,\dots ,n\}$. We denote by $S_n$ the set of permutations on $[n]$ and by $\ST$ the set of surjective maps from $[n]$ to $[r]$. For $n\geq 1$, let ${\mbox {\bf ST}_n}:={\displaystyle \bigcup _{r=1}^n\ST}$.

Note that ${\mbox{\bf ST} _n^n}$ coincides with the set $S_n$ of permutations of $n$ elements, while ${\mbox{\bf ST} _n^1}=\{ c_n\}$, where $c_n$ is the constant function $c_n(i)=1$, for $1\leq i\leq n$. 
\medskip

For $f\in \ST$, we write $\vert f\vert=n$ and  $f = (f(1),\dots ,f(n))$. The composition of maps is denoted $\circ$.

The {\it  concatenation product} $\t : {\mbox {\bf ST}_n^r}\t {\mbox {\bf ST}_m^s}\longrightarrow {\mbox {\bf ST}_{n+m}^{r+s}}$ is given by the formula:
$$f\t g := (f(1),\dots , f(n), g(1)+r,\dots ,g(m) + r).$$
\medskip

\begin{notation}\label{1yepsilon} {\rm We denote by $1_n$ the identity of $S _n$. For any pair of positive integers $n$ and $m$, let $\epsilon (n,m)$ denote the permutation of $n+m$ elements whose image is $(m+1,\dots ,m+n,1,\dots ,m)$. For a finite collection of positive integers $r_1,\dots ,r_s$ with $s>2$, the permutation $\epsilon (r_1,\dots ,r_s)$ in $S_{r_1+\dots +r_s}$ is the composition:
$$\epsilon (r_1,\dots ,r_s):=\epsilon (r_1+\dots +r_{s-1},r_s)\circ(\epsilon (r_1,\dots ,r_{s-1})\times 1_{r_s}).$$
}\end{notation}
\medskip

\begin{definition} \label{std} Given a map $f:[n]\longrightarrow {\mathbb N}$ there exists a unique surjective map ${\mbox {std}(f)}$ in $\ST$ such that $f(i) <f (j)$ if, and only if, ${\mbox {std}(f)}(i) < {\mbox {std}(f)}(j)$, for $1\leq i,j\leq n$. The map ${\mbox {std}(f)}$ is called the {\it standardization} of $f$ (see for instance \cite{Novelli-Thibon}). \end{definition}

For example, when $f=(1,5,4,7,5)$, we get ${\mbox  {std}(f)}=(1,3,2,4,3)$.
\medskip

\begin{notation}\label{notationcorestriction} {\rm For $x\in {\mbox {\bf ST}_n^r}$ and $J=\{ j_1<\dots <j_k\}\subseteq \{1,\dots ,n\}$, let $x\vert_J: ={\mbox {std}(x(j_1),\dots ,x(j_k))}$ denote the restriction of $x$ to $J$. 

Similarly, for $K=\{ j_1<\dots <j_l\}\subseteq \{1,\dots ,r\}$, the co-restriction of $x$ to $K$ is denoted $x\vert ^K:={\mbox {std}(x(s_1),\dots ,x(s_q))}$, for $x^{-1}(K)=\{ s_1<\dots <s_q\}$.

For an element $x\in {\mbox {\bf ST}_n^r}$, we denote by $\lambda (x)$ the cardinal of $x^{-1}(\{r\})$. 

Suppose that $x^{-1}(r)=\{ j_1<\dots <j_{\lambda (x)}\}$, and let $x\rq \in {\mbox {\bf ST}_{n-k}^{r-1}}$ be the co-restriction $x\rq :=x\vert ^{\{1,\dots ,r-1\}}$. We denote $x$ as $x=\prod _{ j_1<\dots <j_{\lambda (x)}}x\rq $.}\end{notation}
\medskip

\begin{example} {\rm For example, the surjective map $x= (3,1,2,5,1,4,3,5,4,2)$ is written as $x=\prod _{4<8}(3,1,2,1,4,3,4,2)$.}\end{example}

\begin{definition}\label{Malpha} Let $x=\prod _{ j_1<\dots <j_{\lambda(x)}}x\rq $ be an element of $\ST$, define the integer ${\mathbb M}(x)_i$, for $0\leq i\leq {\lambda(x)}$ as follows:
$${\mathbb M}(x)_i :=\begin{cases} j_1 - 1,& {\rm for}\ i = 0,\\
j_i - j_{i-1}-1,& {\rm for}\ 1\leq i\leq \lambda(x) -1,\\
n-j_{\lambda (x)},& {\rm for}\ i =  \lambda(x).\end{cases}$$

We define ${\mathbb M}(x):= ({\mathbb M}(x)_{\lambda (x)},\dots , {\mathbb M}(x)_1)$.\end{definition}
\medskip

\begin{definition}\label{xl} {\rm Let $x\in {\mbox {\bf ST}}_n^r$ be a surjective map and let ${\underline l} = (l_1,\dots ,l_p)$ be a collection of integers such that $0=l_0<  l_1 <\dots < l_p < n$.  Define the element $$x^{\underline l} := x\vert ^{\{1,\dots ,l_1\}}\times x\vert ^{\{l_1+1,\dots ,l_2\}}\times \ldots \times x\vert ^{\{l_p+1,\dots ,n\}}.$$ 
 For $p = 0$, define $x^{\underline l} = x$. }\end{definition}
\bigskip

Recall that a {\it composition} of $n$ is a collection $(n_1,\dots ,n_s)$ of positive integers such that $\sum n_i=n$.  

\begin{definition} \label{shuffles} Let $(n_1,\dots ,n_p)$ be a composition of $n$. An element in $f \in {\mbox {\bf ST}_n}$ is a {\it $(n_1,\dots ,n_p)$-stuffle} if 
$$f(n_1+\dots +n_i+1) < f(n_1+\dots +n_i+1) < \dots < f(n_1+\dots +n_i+n_{i+1}),$$
for $0\leq i\leq p-1$. \end{definition}

\begin{notation}\label{notation3} We denote by ${\mbox {\it SH}(n_1,\dots ,n_p)}$ the set of all $(n_1,\dots ,n_p)$-stuffles.
\medskip

For a composition $(n_1, \dots , n_p)$ of $n$, we denote: \begin{enumerate}
\item  ${\mbox {SH}^{\prec}(n_1,\dots ,n_p)}$ the subset of all surjective maps $f\in {\mbox {SH}(n_1,\dots ,n_p)}$ such that $f(n_1)>f(n_1+n_2)>\dots >f(n)$.
\item  ${\mbox {SH}^{\succ}(n_1,\dots ,n_p)}$ the subset of all surjective maps $f\in {\mbox {SH}(n_1,\dots ,n_p)}$ such that $f(n_1)<f(n_1+n_2)<\dots <f(n)$.
\item ${\mbox {SH}^{\bullet}(n_1,\dots ,n_p)}$ the subset of all surjective maps $f\in {\mbox {SH}(n_1,\dots ,n_p)}$ such that $f(n_1)=f(n_1+n_2)=\dots =f(n)$.
\item ${\mbox {SH}^{\succcurlyeq}(n_1,\dots ,n_p)}$ the subset of all surjective maps $f\in {\mbox {SH}(n_1,\dots ,n_p)}$ such that $f(n_1)\leq f(n_1+_2)\leq \dots \leq f(n)$.\end{enumerate}
\end{notation}
\medskip

To recover the usual notion of shuffle, it suffices to note that a $(n_1,\dots ,n_p)$-shuffle is a permutation $\sigma \in S_n\cap {\mbox {\it SH}(n_1,\dots ,n_p)}$. We denote by ${\mbox {\it Sh}(n_1,\dots ,n_p)}$ the set of all $(n_1,\dots ,n_p)$-shuffles.

In an analogous way, we define ${\mbox {Sh}^{\prec}(n_1,\dots ,n_p)}: = S_n \cap  {\mbox {SH}^{\prec}(n_1,\dots , n_p)}$, and ${\mbox {Sh}^{\succ}(n_1,\dots ,n_p)}:= S_n \cap  {\mbox {SH}^{\succ}(n_1,\dots , n_p)}$.

Finally, we denote by ${\mbox {Sh}^{\bullet}(r_1,\dots ,r_p)}$ the set of all $f = {\displaystyle \prod _{r_1 <r_1+r_2<\dots < r_1+\dots +r_p}f\rq }$ with $f\rq \in {\mbox {Sh}(r_1{-}1,\dots ,r_p {-}1)}$.
\medskip

The following property of the shuffles is well-known and is the key result to prove the associativity of the shuffle product.

\begin{proposition} \label{shuffleassociativity} Let $n, m$ and $r$ be positive integers. The set of $(n,m,r)$-shuffles satisfies the following property:
$${\mbox {Sh}(n+m,r)}\circ ({\mbox {Sh}(n,m)}\times 1_r) = {\mbox {Sh}(n,m,r)} = {\mbox {Sh}(n, m+r)}\circ (1_n\t {\mbox {Sh}(m,r)}),$$
where $1_n= (1,2,\dots ,n)$ denotes the identity of $S_n$. \end{proposition}
\medskip

\begin{remark} \label{shuffledendr} {\rm A straightforward calculation shows that the equality of Proposition \ref{shuffleassociativity} splits into three formulas:\begin{enumerate}
\item ${\mbox {Sh}^{\succ}(n+m,r)}\circ ({\mbox {Sh}(n,m)}\times 1_r) = {\mbox {Sh}^{\succ}(n, m+r)}\circ (1_n\t {\mbox {Sh}^{\succ}(m,r)}),$
\item ${\mbox {Sh}^{\prec}(n+m,r)}\circ ({\mbox {Sh}^{\succ}(n,m)}\times 1_r) = {\mbox {Sh}^{\succ}(n, m+r)}\circ (1_n\t {\mbox {Sh}^{\prec}(m,r)}),$
\item ${\mbox {Sh}^{\prec}(n+m,r)}\circ ({\mbox {Sh}^{\prec}(n,m)}\times 1_r) = {\mbox {Sh}^{\prec}(n, m+r)}\circ (1_n\t {\mbox {Sh}(m,r)}).$\end{enumerate}}\end{remark}

The set of stuffles ${\mbox {SH}(n,m,r)}$ satisfies analogous properties. We shall use them to define tridendriform algebra structures, for details we refer to \cite{Palacios-Ronco}.

\bigskip

For $n\geq 1$ and $1\leq i\leq n-1$, let $t_i\in S_n$ be the permutation which exchanges $i$ and $i+1$, that is 
$$t_i(j) := (1,\dots ,i-1,i+1,i,i+2,\dots ,n) $$

\begin{definition} \label{wBo} For $n\geq 1$, the weak Bruhat order on the set $S_n$ of permutations is defined by the covering relation:
$$\sigma < t_i\circ \sigma,$$
when $\sigma ^{-1}(i) < \sigma ^{-1}(i+1)$.\end{definition} 

The following Proposition is well-known, see for instance \cite{Aguiar-Sottile} or \cite{Loday-Roncoord}.

\begin{proposition}  \label{orderST} For any composition $(n_1,\dots ,n_p)$ of $n$, the set of shuffles ${\mbox {Sh}(n_1,\dots ,n_p)}$ coincides with the subset $\{ w\in S_n\mid 1_n\leq w \leq \epsilon (n_1,\dots ,n_p)\}$ of $S_n$, where $\leq $ is the weak Bruhat order. \end{proposition}

The weak Bruhat order of $S_n$ may be extended to the set of surjective maps ${\mbox {\bf ST}_n}$ in two different ways (see \cite{Palacios-Ronco}). We describe the one we need in the last section of the paper.

\begin{definition}\label{PalRo} For $n\geq 1$, the weak Bruhat order on ${\mbox {\bf ST}_n^r}$ is the transitive relation spanned by the covering relation $$f < t_i\circ f,\ 
{\rm when}\ f ^{-1}(i) <f^{-1}(i+1),$$
for some $1\leq i\leq r-1$, where for any pair of subsets $J,K\subseteq \{1,\dots ,r\}$ we say that $J < K$ if the maximal element of $J$ is smaller that the minimal element of $K$.
\end{definition}

For instance, $(1, 4 , 1, 3, 4, 2) < (2,4,2,3,4,1)$ , but the elements $(1, 4 , 1, 3, 4, 2)$ and $(1, 3 , 1, 4, 3, 2)$ are not comparable.
\medskip

The following result is proved in \cite{Palacios-Ronco}.

\begin{proposition} \label{weakBST} Let $\sigma < \tau\in {\mbox {Sh}(r_1,\dots ,r_p)}$ be two permutations, and let $x_i \leq  x_i\rq \in {\mbox {\bf ST}_{n_i}^{r_i}}$, for $1\leq i\leq p-1$ be surjective maps. We have that:\begin{enumerate}
\item $\sigma \circ (x_1\times \ldots \times x_p) < \tau  \circ (x_1\times \ldots \times x_p)$,
\item $\sigma \circ (x_1\times \ldots \times x_p) \leq \sigma \circ (x_1\rq \times \ldots \times x_p\rq )$. Moreover, if at least for one $1\leq i\leq p$ the elements $x_i$ and $x_i\rq$ satisfy that $x_i < x_i\rq$, then $\sigma \circ (x_1\times \ldots \times x_p)< \sigma \circ (x_1\rq \times \ldots \times x_p\rq )$.\end{enumerate}\end{proposition}

\bigskip

\bigskip

\section{The coalgebra $\st$ of surjective maps}
\medskip

We define different algebraic structures on the graded vector space $\K[{\mbox {\bf ST}}]:=\bigoplus _{n\geq 1}\K[{\mbox {\bf ST}_n}]$. We begin by the coassociative coproduct $\Delta $. For a more detailed description of the properties of $\Delta $ see \cite{Novelli-Thibon} , \cite{Chapoton} or \cite{Loday-Ronco}. 
\medskip

\begin{definition} \label{coproduct} {\rm We define $\Delta : \K[{\mbox {\bf ST}}]\longrightarrow \K[{\mbox {\bf ST}}] \ot \K[{\mbox {\bf ST}}]$ on an element $x\in \ST$ by:
$$\Delta (x) = \sum _{i=1}^{r-1} x\vert ^{\{1,\dots , i\}}\ot x\vert ^{\{ i+1, \dots , r\}},$$
and we extend it by linearity to all $\K[{\mbox {\bf ST}}]$}\end{definition}

For example,
$$\displaylines{
\Delta (3, 4, 2, 5, 1, 1, 3, 5) = (1, 1)\ot (2, 3, 1, 4, 2, 4) + (2, 1, 1) \ot (1, 2, 3, 1, 3) +\hfill\cr
\hfill (3, 2, 1, 1, 3)\ot (1, 2, 2) + (3, 4, 2, 1, 1, 3)\ot (1, 1).\cr }$$

For any $x\in \ST $ and any pair $1< i < j < r-1$, we have that:\begin{enumerate}
\item $(x\vert ^{\{1, \dots , j\}})\vert ^{\{i+1,\dots ,j\}} = x\vert ^{\{ i+1, \dots , j\}}$,
\item $(x\vert ^{\{ i+1, \dots , r\}})\vert ^{\{1,\dots ,j\}} =  x\vert ^{\{ i+1, \dots , j\}}$,\end{enumerate}
which implies that the coproduct is coassociative.
\medskip 

Let $\st$ denote the  graded coalgebra $(\K[{\mbox {\bf ST}}], \Delta)$. 
On $\st ^+$, the coproduct $\Delta $ is uniquely extended to $\Delta ^+$ in such a way that the reduced coproduct of $\st ^+$ is $\Delta$.

The data $(\st ^{+}, \Delta ^{+})$ is a coassociative unital and counital coalgebra.
\bigskip

It is clear that the concatenation product $\t$, extended by linearity, defines an associative graded product on $\st$.

\begin{definition}\label{irreducibles} {\rm An element $f\in \ST$ is called {\it irreducible} if there do not exist an integer $1\leq i\leq n-1$ and a pair of surjective maps $g\in{\mbox {\bf ST}_{i}^{k}}$ and $h \in{\mbox {\bf ST}_{n-i}^{r-k}}$ such that $f = g\t h$. We denote by ${\mbox {\bf Irr}_n}$ the set of irreducible elements of ${\mbox {\bf ST}_n}$, for $n\geq 1$, and by ${\mbox {\bf Irr}}$ the union $\bigcup _{n\geq 1}{\mbox {\bf Irr}_n}$.}\end{definition}
\medskip

\begin{remark} \label{irreduciblesurj} Given a surjective map $x\in {\mbox {\bf ST}_n^r}$, there exists a unique family $x^1,\dots ,x^p$ of elements, with $x^i\in {\mbox {\bf Irr}_{n_i}^{r_i}}$, such that 
$x = x^1\times \ldots \times x^p$, where $n=\sum _{i=1}^pn_i$ and $r=\sum _{i=1}^pr_i$. So, the space $\K[{\mbox {\bf ST}}]$ with  $\times $ is the free associative algebra spanned by the set ${\mbox {\bf Irr}}$.\end{remark}
\medskip

\begin{lemma} \label{lemma0} Let $x$ and $y$ be elements of $\st ^+$, we have that:
$$\Delta^+ (x\t y) = \sum \bigl( x_{(1)}\ot (x_{(2)}\t y) + (x\t y_{(1)})\ot y_{(2)}\bigr ) - x\ot y,$$
where $\Delta^+(x) = \sum x_{(1)}\ot x_{(2)}$, $\Delta ^+ (y) = y_{(1)}\ot y_{(2)}$.\end{lemma}
\medskip

\begin{proo} Suppose that $x\in {\mbox {\bf ST}_n^r}$ and $y\in {\mbox {\bf ST}_m^s}$. To prove the Lemma it suffices to note that: \begin{enumerate}
\item $$(x\t y)\vert ^{\{ 1, \dots ,i\}} = \begin{cases} x\vert^{\{1,\dots ,i\}}, & {\rm for}\ 0\leq i\leq r,\\
x\t y\vert^{\{ 1,\dots , i - r \}}, & {\rm for}\ r+1\leq i\leq r+s,\end{cases}$$
\item $$(x\t y)\vert ^{\{ i+1, \dots ,r+s\}} = \begin{cases} x\vert^{\{i+1,\dots ,r\}} \t y, &{\rm for}\ 0\leq i\leq r,\\
y\vert ^{\{ i-r+1,\dots , s\}}, &{\rm for}\ r+1\leq i\leq r+s,\end{cases}$$
\end{enumerate}
which imply the formula.\end{proo}
\bigskip

A vector space $V$ equipped with an associative product and a coassociative coproduct, satisfying the condition of Lemma \ref{lemma0} is called an {\it infinitesimal unital bialgebra} in \cite{Loday-Ronco2}. So, $(\st^+, \t , \Delta^+)$ is a unital infinitesimal bialgebra. 
\medskip

As proved in \cite{Loday-Ronco2}, any conilpotent unital infinitesimal bialgebra $(C, \t, \Delta )$ is isomorphic, as a coalgebra, to the cotensor algebra $T^c({\mbox {Prim}(C)})$. Moreover, the linear map 
$$E (x) := \sum _{i\geq 1} (-1)^i \bigl (\sum x_{(1)} \t \ldots \t x_{(i)}\bigr ),$$ gives a projection from $C$ to ${\mbox {Prim}(C)}$, where 
${\overline {\Delta}}^i(x) = \sum x_{(1)} \ot \dots \ot x_{(i)}$ for $x\in C$. For the details of the construction we refer to  \cite{Loday-Ronco2}. 

\begin{remark} \label{idempotentE} For the particular case of $(\st ^+, \t , \Delta^+)$, we get that the linear map $E: \K[{\mbox {\bf Irr}}] \longrightarrow {\mbox {Prim}(\st^+)}$ is an isomorphism. It induces an isomorphism of coalgebras $E_{ST}: T^c(\K[{\mbox {\bf Irr}}]) \longrightarrow \st$, given by:
$$E_{ST}(x^1\ot \dots \ot x^p):= E(x^1)\t \dots \t E(x^p),$$
for any family of irreducible elements $x^1,\dots ,x^p\in {\mbox {\bf Irr}}$. On the other hand, we proved that $Id_{\st_n} = \sum _{j=1}^n\t ^j\circ E^{\ot j}\circ {\overline {\Delta}}^j$, which implies that $T^c({\mbox {Prim}(\st)})$ is isomorphic to $\st$, as coalgebras.  Thus, putting the isos together, we have
$$T^c(\K[{\mbox {\bf Irr}}])\cong \st \cong T^c({\mbox {Prim}(\st )})\;.$$
Noting that these isos respect the grading induced on  the tensor algebras, we get 
that  the dimension of the subspace $\mbox{Prim}(\st)_n
$ of homogeneous elements of degree $n$ is the cardinal $\vert {\mbox {\bf Irr}_n}\vert $  of the set of irreducible elements of degree $n$, for $n\geq 1$.\end{remark}
\medskip

Using Definition \ref{xl}, $E(x)= \sum_{\underline l} \alpha _{\underline l}x^{\underline l}$, where the sum is taken over all families ${\underline l} = (l_1,\dots ,l_p)$ such that $0=l_0< l_1 < \ldots < l_p < n$, $p\geq 0$ and $\alpha _{\underline l}:= (-1)^p$, which implies that $x^{\underline l} = x$ or $x^{\underline l}$ is reducible. 
\medskip

\begin{remark} \label{partition} {\rm If $x\in {\mbox {\bf ST}_n^r}$ and ${\underline l} = (l_1,\dots ,l_p)$, then 
$x^{\underline l}\in {\mbox {\bf ST}_n^r}$, too.}\end{remark}
\medskip

\begin{lemma}\label{lemmaMalphaE} Let $x\in {\mbox {\bf ST}_n^r}$ and let $p\geq 1$. Given a family of integers ${\underline l}$ such that $0 = l_0 < l_1 < l_2< \ldots < l_p < r$, the element 
$x^{\underline l}$ satisfies that ${\mathbb M}(x^{\underline l})\leq {\mathbb M}(x)$ for the lexicographic order. \end{lemma}
\medskip 

\begin{proo} Suppose that $x = \prod _{j_1<\ldots <j_{\lambda (x)}}x\rq$. 

If $x^{-1}(\{1,\dots ,l_p\})\subseteq \{1,\dots ,j_1-1\}$, then $x^{\underline l} = \prod _{j_1< \ldots < j_{\lambda (x)}}(x\rq )^{\underline l}$, which implies that  ${\mathbb M}(x^{\underline l}) = {\mathbb M}(x)$.

On the other hand, if  $x^{-1}(\{1,\dots ,l_p\})\cap \{j_1+1,\dots ,n\}\neq \emptyset $, then there exists at least one $j_1< k \leq n$ such that $x(k)\leq l_p$. Let $k_0$ be the maximal integer which satisfies this condition. There exists $1\leq i_0 \leq p$ such that $j_{i_0} < k_0 <j_{i_0 +1}$, and we get that: \begin{itemize}
\item ${\mathbb M}(x^{\underline l}) _i = {\mathbb M}(x)_i$ for $i_0 < i \leq \lambda (x)$,
\item ${\mathbb M}(x^{\underline l}) _{i_0} < {\mathbb M}(x)_{i_0}.$\end{itemize}
So, ${\mathbb M}(x^{\underline l} ) < {\mathbb M}(x)$.
\end{proo}

\bigskip

\bigskip

\section{The dendriform bialgebra  $\st_D$}\label{dendbialg}
\medskip

The shuffle product defines a bialgebra structure on $\st ^+$, which has been studied by F. Chapoton, and by J.-C. Novelli and J.-Y. Thibon, who called it the bialgebra of {\it packed words}. We recall the main constructions and results, for the details of the proofs we refer to \cite{Chapoton} and \cite{Novelli-Thibon}.

Let $x\in {\mbox {\bf ST}_n^r}$ and $y\in {\mbox {\bf ST}_m^s}$ be two surjective applications, the {\it shuffle product}  $x * y$ is defined by:
$$x*y:=\sum _{f\in Sh(r,s)} f\circ (x\t y).$$

The product $*$ is associative and satisfies that:
$$\Delta (x*y) = \sum (x_{(1)}*y_{(1)})\otimes (x_{(2)}*y_{(2)}),$$
for any $x, y \in {\mbox {\bf ST}}$. So, the coalgebra $(\st^+, \Delta )$ equipped with the shuffle product $*$ is a bialgebra over $\K$.
\medskip

Using Remark \ref{shuffledendr}, the shuffle product of $\st^+$ comes from a dendriform structure of $\st$.

That is, the vector space $\st$ with the products $\succ $ and $\prec$ defined by:\begin{enumerate}
\item $x\succ y := \sum _{f\in Sh^{\succ}(r,s)} f\circ (x\t y)$,
\item $x\prec y := \sum _{f\in Sh^{\prec}(r,s)} f\circ (x\t y)$,\end{enumerate}
for $x\in \ST$ and $y\in {\mbox {\bf ST}_m^s}$, is a dendriform algebra.
\medskip

Fix that $x\succ 1_{\K} := 0 =: 1_{\K}\prec x$ and $x\prec 1_{\K} := x =: 1_{\K}\succ x$, for all $x\in \st$. Note that $\st^+$ is not a dendriform algebra, because there does not exist a coherent way to define $1_{\K}\succ 1_{\K}$ and $1_{\K}\prec 1_{\K}$. It is easily seen that $\st $ is a conilpotent dendriform bialgebra. 
\medskip

The dendriform algebra $(\st , \succ ,\prec )$ is free. The main result of the present section is to give a proof of this result by exhibiting a basis.

\begin{notation} \label{notationunder} {\rm Let $x\in {\mbox {\bf ST}_n^r}$ and $y\in {\mbox {\bf ST}_m^s}$ be two maps, we denote by $x\backslash y$ the composition:
$$x\backslash y : = \epsilon (r,s)\circ (x\times y)\in {\mbox {\bf ST}_{n+m}^{r+s}}.$$}\end{notation}
\medskip

For $n\geq 1$, define the subset ${\mathcal D}_n\subseteq {\mbox {\bf ST}_n}$ recursively, as follows:\begin{enumerate}
\item ${\mathcal D}_1 := \{ (1)\} = {\mbox {\bf ST}_1}$,
\item ${\mathcal D}_2 := \{ (1,1)\} $,
\item a surjective map $x\in {\mbox {\bf ST}_n}$ belongs to ${\mathcal D}_n$ if $x$ is irreducible, and there do not exist an integer $1\leq r < n$ and a pair of elements $y\in {\mathcal D}_r$ and $z\in {\mbox {\bf ST}_{n-r}}$ such that $x = y\backslash z$.\end{enumerate}

\begin{theorem}\label{baseasdend} The dendriform algebra $(\st ,\succ , \prec)$ is the free dendriform algebra spanned by the set 
$$E({\mathcal D}) = \{ E(x)\mid x\in \bigcup _{n\geq 1}{\mathcal D}_n\},$$
 where $E(x) =  \sum _{i\geq 1} (-1)^i \bigl (\sum x_{(1)} \t\dots \t x_{(i)}\bigr )$, with ${\overline \Delta} ^i(x) = \sum x_{(1)}\ot \dots \ot x_{(i)}$, for $i\geq 1$.\end{theorem}
\medskip

In order to prove Theorem \ref{baseasdend}, we need some additional results. Note first that:

\begin{remark} \label{irredandD} {\rm  For any irreducible element $x\in{\mbox {\bf Irr}}$, there exist a unique integer $1\leq m\leq n$ and unique  elements $y\in  {\mathcal D}_m$ and $z\in {\mbox {\bf ST}_{n-m}}$ such that $x = y\backslash z$. Moreover, there exists a unique way to write down $z = z^1\t \dots \t z^p$, with $z^i\in {\mbox {\bf Irr}}$, for $1\leq i\leq p$.}\end{remark}
\medskip

The next Lemmas will serve in  the proof of Theorem \ref{baseasdend}.

\begin{lemma}\label{lemma1} Let $x\in {\mbox {\bf ST}_n^r}$ and $y\in {\mbox {\bf ST}_m^s}$. For any $f\in {\mbox {Sh}(r,s)}$, we have that:\begin{enumerate}
\item $f\circ (x\times y)\vert _{\{1,\dots ,n\}}= x$,
\item $f\circ (x\times y)\vert _{\{n+1,\dots ,n+m\}} = y$.\end{enumerate}\end{lemma}
\medskip

\begin{proo} The result is an easy consequence of 
$$f(1)<\dots <f(r),\ {\rm and}\ f(r+1)<\dots < f(r+s).$$\end{proo}
\medskip

\begin{notation}\label{notationbracel} {\rm We denote by ${\mathfrak B}(l)$ the set of all the elements $x$ in ${\mbox {\bf ST}}$ which are of the form $x = y\backslash z$ with $y\in {\mathcal D}$ and $\vert z \vert =l$, for $l\geq 0$.  Note that ${\mathfrak B}(0) = {\mathcal D}$.} \end{notation}
\medskip

\begin{lemma}\label{lemma2} Let $x\in {\mathcal D}_n^r$ and $y\in {\mbox {\bf ST}_m^s}$ be two surjective maps and let $f\in {\mbox {Sh}^{\prec}(r,s)}$, $f\neq \epsilon (r,s)$. If there exist $z\in {\mathcal D}$ and $w\in\mbox {\bf ST}$  such that $f\circ (x\t y) = z\backslash w\in {\mathfrak B}(l)$, then $l < m$.\end{lemma}
\medskip

\begin{proo}  Suppose that $f\circ (x\t y) = z\backslash w$ with $\vert w\vert  \geq m$. In this case, $\vert z\vert \leq n$, and from Lemma \ref{lemma1} we get that:
$$x = f\circ (x \times y)\vert _{\{1,\dots , n\}}  = z\backslash (w\vert _{\{1,\dots , n-\vert z\vert\}}).$$

As $x\in {\mathcal D}$, then $x = z$ and $n =\vert z\vert$. So, we get $f\circ (x\times y ) = x\backslash w$, but this is possible only when $f = \epsilon (r,s)$.\end{proo}
\medskip

\begin{lemma} \label{lemma3} Let $x$ be a reducible element of ${\mbox {\bf ST}_n^r}$ and let $y\in {\mbox {\bf ST}_m^s}$. For any $f\in {\mbox {Sh}^{\prec}(r,s)}$, we have that either 
$f\circ (x\times y)$ is reducible, or $f\circ (x\times y) \in {\mathfrak B}(l)$, for $l < m$.\end{lemma}
\medskip

\begin{proo} Suppose that $f\circ (x\times y)$ is irreducible. In this case, there exist $z\in {\mathcal D}$ and $w\in {\mbox {\bf ST}}$ such that $f\circ (x\times y) = z\backslash w$.

Lemma \ref{lemma1} states that $x = f\circ (x\times y)\vert _{\{1,\dots ,n\}} =(z\backslash w)\vert_{\{1,\dots ,n\}}$. If $\vert z\vert \leq n$, we get that $z$ is reducible, which is false. So, $\vert z\vert  > n$, which implies that $f\circ (x\times y) = z\backslash w \in {\mathfrak B}(l)$, for $l < m$.\end{proo}
\bigskip

Define, on the vector space $\K[{\mbox {\bf Irr}}]$, a structure of brace algebra given by:\begin{enumerate}
\item for $x\in {\mathcal D}$ and $y_1,\dots ,y_n \in {\mbox {\bf Irr}}$,
$$M_{1n}(x; y_1,\dots ,y_n) := x\backslash (y_1\t \dots \t y_n),$$
\item for $x\in {\mbox {\bf Irr}}\setminus {\mathcal D}$ and and $y_1,\dots ,y_n \in {\mbox {\bf Irr}}$,  there exist $x_1\in {\mathcal D}$ and $x_2\in {\mbox {\bf ST}}$ such that $x_2 = z_1\t \ldots \t z_p$ with $z_j\in {\mbox {\bf Irr}}$, for $1\leq j\leq p$. 

\noindent In this case, $M_{1n}(x; y_1,\dots ,y_n) = M_{1n}(M_{1p}(x_1; z_1,\dots, z_p);y_1,\dots ,y_n)$ is defined using Definition \ref{def:Gerstenhaber-Voronov}.\end{enumerate}
\medskip
 
Lemma \ref{bracebasis} states that $(\K[{\mbox {\bf Irr}}], \{ M_{1n}\}_{n\geq 1})$ is a well-defined brace algebra. 

A recursive argument on $\vert y\vert $, for $y\in {\mbox {\bf Irr}}$, shows that as a brace algebra $\K[{\mbox {\bf Irr}}]$ is freely generated by ${\mathcal D}$. 
\medskip

Let $E\vert _{ {\mathcal D}}:  {\mathcal D}\longrightarrow {\mbox {Prim}(\st)}$ be the restriction to $ {\mathcal D}$ of the projection $E(x) =  \sum _{i\geq 1} (-1)^i \bigl (\sum x_{(1)} \t\dots \t x_{(i)}\bigr )$.  As ${\mbox {Prim}(\st)}$ is a brace algebra, there exists a unique homomorphism of brace algebras $\eta: \K[{\mbox {\bf Irr}}]\longrightarrow {\mbox {Prim}(\st)}$, such that $\eta (x) = E(x)$, for $x \in {\mathcal D}$. 
 \medskip

We want to prove that $\eta $ is an isomorphism.  
\medskip

\begin{notation}\label{notationeta} {\rm Let $x\in {\mbox {\bf Irr}}$, we denote  $\eta (x) = \sum_{a_i\neq 0} a_ix_{i}$.  There exists a unique $i_0$ such that $x_{i_0}= x$ and $a_{i_0}=1$.} \end{notation}
\medskip

Given an element $x = y\backslash z\in{\mbox {\bf Irr}_n}$, with $y\in  {\mathcal D}_m$ and $z = z^1\t \dots \t z^p$, with $z^i\in {\mbox {\bf Irr}}$ for $1\leq i\leq p$, the homomorphism $\eta$ is defined as 
$$\displaylines {
\eta (x) = M_{1p}^{ST}(\eta (y); \eta (z^1),\dots ,\eta (z^p)) =\hfill\cr
\hfill \sum _j(-1)^j\bigl (\sum  c_{k,i_1,\dots ,i_p}\omega ^{\prec}(z_{i_1}^1,\dots ,z_{i_j}^j)\succ y^{\underline l}\prec \omega ^{\succ}(z_{i_{j+1}}^{j+1},\dots z_{i_p}^p)),\cr} $$
where $M_{1p}^{ST}$ denote the brace operations in $\st$ and $\eta (z^j) =\sum _{b_{i_j}\neq 0} b_{i_j} z_{i_j}^j$, for $1\leq j\leq p$. 
\medskip

\begin{lemma} \label{lemma4} Let $x = y\backslash z$ be an element in ${\mathfrak B}(l)$, with $y\in {\mathcal D}_n^r$. 

\noindent If $\eta (x) = \sum _{a_i\neq 0}a_i x_{i}$, then:\begin{enumerate}
\item ${\mathbb M}(x_{i}) \leq {\mathbb M}(x)$ for the lexicographic order,
\item if ${\mathbb M}(x_{i}) = {\mathbb M}(x)$, then $x_{i}$ is reducible or $x_{i}\in {\mathfrak B}(k)$, with $0\leq k\leq l$,\end{enumerate}
for all $i$. 

Suppose that $z = z^1\t \dots \t z^p$, with $z^j\in {\mbox {\bf Irr}_{m_j}^{s_j}}$ and $\eta (z^j) = \sum _{b_{i_j}\neq 0}b_{i_j}z_{i_j}^j$. The unique terms $x_{i}$ satisfying that ${\mathbb M}(x_{i}) = {\mathbb M}(x)$ and $x_{i}\in {\mathfrak B}(l)$, are of the form:
$$x_{i} = y\backslash (g\circ (z_{i_1}^{1},\dots , z_{i_p}^p)),$$
with $g\in {\mbox {Sh}^{\succ}(s_{1},\dots ,s_p)}$, for some family $z_{i_j}^{j}$.\end{lemma}
\medskip

\begin{proo} Note first that $\lambda (x) = \lambda (y)$ and ${\mathbb M}(x)_{\lambda (x)} = {\mathbb M}(y)_{\lambda (y)} + \vert z\vert $.

As $y\in {\mathcal D}$, we have that $\eta (y) = \sum _{\underline l} \alpha_{\underline l} y^{\underline l}$, with ${\mathbb M}(y^{\underline l})\leq {\mathbb M}(y)$ and $y^{\underline l}$ reducible or $y^{\underline l}=y$, by Lemma \ref{lemmaMalphaE}.
\medskip

Recall that 
$$\displaylines { 
\omega ^{\prec}(z_{i_1}^1,\dots , z_{i_j}^j)\succ y^{\underline l}\prec \omega ^{\succ}(z_{i_{j+1}}^{j+1},\dots , z_{i_p}^p)=\hfill\cr
\sum _{f, g, h} h\circ (f\circ (z_{i_1}^1\t \dots \t z_{i_j}^j)\t y^{\underline l} \t g\circ (z_{i_{j+1}}^{j+1}\t \dots \t z_{i_p}^p),\cr }$$
where  $f\in {\mbox {Sh}^{\prec}(s_1,\dots ,s_j)}$, $g\in {\mbox {Sh}^{\succ}(s_{j+1},\dots ,s_p)}$ and 

\noindent $h\in {\mbox {Sh}(s_1+\dots +s_j, r,s_{j+1}+\dots + s_p)}$ is such that $h(s_1+\dots +s_j) < h(s_1+\dots +s_j + r)$ and $h(s_1+\dots +s_j + r) > h(s_1+\dots +s_p+ r)$.

Let $x_{i} := h\circ (f\circ (z_{i_1}^1\t \dots \t z_{i_j}^j)\t y^{\underline l} \t g\circ (z_{i_{j+1}}^{j+1}\t \dots \t z_{i_p}^p))$.

\noindent We have that $\lambda (x_{i}) = \lambda (y^{\underline l})$ and ${\mathbb M}(x_{i})_{\lambda (x_{i})} = {\mathbb M}(y^{\underline l})_{\lambda (y^{\underline l})} + \vert z^{j+1}\vert + \dots + \vert z^p\vert $. So, \begin{enumerate}
\item ${\mathbb M}(x_{i}) \leq {\mathbb M}(x)$, for all $x_{i}$,
\item if $j\geq 1$, then ${\mathbb M}(x_{i}) < {\mathbb M}(x)$,
\item if ${\mathbb M}(y^{\underline l}) < {\mathbb M}(y)$, then ${\mathbb M}(x_{i}) < {\mathbb M}(x)$.\end{enumerate}

The unique elements such that ${\mathbb M}(x_{i}) = {\mathbb M}(x)$ are of the form $$x_{i} = h\circ (y^{\underline l}\t g\circ (z_{i_1}^{1},\dots , z_{i_p}^p)),$$ where
${\mathbb M}(y^{\underline l}) = {\mathbb M}(y)$, $g\in {\mbox {Sh}^{\succ}(s_{1},\dots ,s_p)}$ and $h\in {\mbox {Sh}^{\prec}(r, s_1+\dots +s_p)}$.
\medskip

But, if $y^{\underline l}\neq y$, then $y^{\underline l}$ is reducible, and by Lemma \ref{lemma3}, we have that either $x_{i}$ is reducible or $x_{i}\in {\mathfrak B}(k)$, with $0\leq k< l$. So, we may restrict ourselves to consider only the $x_{i}$ of the form:
$$x_{i} = h\circ (y\t g\circ (z_{i_1}^{1},\dots , z_{i_p}^p)),$$
where $g\in {\mbox {Sh}^{\succ}(s_{1},\dots ,s_p)}$ and $h\in {\mbox {Sh}^{\prec}(r, s_1+\dots +s_p)}$.

Lemma \ref{lemma2} implies that for $h\neq \epsilon(r, s_1+\dots +s_p)$, the element  $x_{i}$ belongs to ${\mathfrak B}(k)$, for $k < l$. 

The unique terms such that ${\mathbb M}(x_{i}) = {\mathbb M}(x)$ and $x_{i}\in {\mathfrak B}(l)$, are of the form:
$$x_{i} = y\backslash (g\circ (z_{i_1}^{1},\dots , z_{i_p}^p)),$$
where $g\in {\mbox {Sh}^{\succ}(s_{1},\dots ,s_p)}$.\end{proo}
\medskip

\begin{definition} \label{etaover} {\rm For $x \in {\mbox {\bf Irr}}$, define ${\overline {\eta}}(x)$ as follows:\begin{enumerate}
\item ${\overline {\eta}}(x) := x$, for $x\in {\mathcal D}$,
\item ${\overline {\eta}}(x) := y\backslash \omega^{\succ}(\eta (z^1),\dots ,\eta (z^p))$, 

\noindent for $x = y\backslash z\in {\mathfrak B}(l)$, with $y\in {\mathcal D}$ and $z = z^1\t \dots \t z^p$, such that $z^j\in {\mbox {\bf Irr}}$.\end{enumerate}}\end{definition}
\medskip

\begin{corollary} \label{corollarybaseasdend} If the set $\{{\overline {\eta}} (x)\mid x\in {\mbox {\bf Irr}}\}$ is linearly independent in $\st$, then the set $\{ \eta (x)\mid x\in {\mbox {\bf Irr}}\}$ is linearly independent in ${\mbox {Prim}(\st)}$.\end{corollary}
\bigskip

\centerline{\bf Proof of Theorem \ref{baseasdend}}
\medskip

Recall from  Remark \ref{idempotentE}  that ${\mbox {dim}_{\K}(\K[{\mbox {\bf Irr}}_n])}$ and ${\mbox {dim}_{\K}({\mbox {Prim}(\st )}_n)}$ are equal. Hence it suffices to verify that $\eta$ is either injective or surjective, grade by grade.
\medskip

For any irreducible element $x$, we have that $E(x) = \sum _{{\underline l}} \alpha _{\underline l} x^{\underline l}$, where $x^{\underline l}$ is reducible for all ${\underline l}= (l_1,\dots ,l_p)$, $ p\geq 1$. So, the set $\{ \eta(x)\mid x\in {\mathcal D}\}$ is linearly independent in ${\mbox {Prim}(\st)}$. 
\medskip

Applying Corollary \ref{corollarybaseasdend}, it suffices to show that the set $\{{\overline {\eta}} (x)\mid x\in {\mbox {\bf Irr}}\}$ is linearly independent in $\st$.
\bigskip

As ${\mbox {Prim}(\st)} = \bigoplus _{n\geq 1}{\mbox {Prim}(\st)_n}$, we prove the result by induction on $n$. 

For $n =1$, ${\mbox {Prim}(\st)_1}=\K\cdot (1)$ and $(1) = \eta (1) = {\overline {\eta}}(1)$.
\medskip

For $n = 2$, ${\mbox {\bf Irr}_2} = \{ (2, 1) ; (1, 1)\}$. We have that ${\overline {\eta}} ((2, 1)) = (2, 1)$ and ${\overline {\eta}} (1,1) = (1,1)$, which proves the result.
\medskip

For $n \geq 3$, suppose that ${\overline {\eta}}({\mbox {\bf Irr}_m})$ is linearly independent for all $m < n$. We have 
that ${\eta}({\mbox {\bf Irr}_m})$ is linearly independent in ${\mbox {Prim}(\st)}$ for all $m < n$, which implies that ${\displaystyle \bigoplus _{j=1}^{n-1}{\mbox {Prim}(\st)_j}}$ is spanned 
by ${\displaystyle \bigcup _{j=1}^{n-1}\eta({\mbox {\bf Irr}_j})}$. 

Consider the linear map ${\overline {\varphi}}: \st\longrightarrow \st$ defined by:
$${\overline {\varphi}}(z) := \omega ^{\succ}(\eta (z^1),\dots ,\eta(z^p)),$$
for $z = z^1\t\dots \t z^p$ with $z^1,\dots ,z^p$ irreducible elements. 
\medskip

Proposition \ref{propositionbrace} asserts that any element of ${\displaystyle \bigoplus _{j=1}^{n-1}\st _j}$ belongs to 

$${\overline {\varphi}}(\bigoplus _{j=1}^{n-1}{\mbox {Prim}(\st)_j}) = {\overline {\varphi}}(\K[\bigcup _{j=1}^{n-1}\eta({\mbox {\bf Irr}_j)}]).$$

Therefore the set $\{ {\overline {\varphi}}(z)\ \mid z\in {\displaystyle \bigcup _{j=1}^{n-1}{\mbox {\bf ST}_j}}\}$ spans ${\displaystyle \bigoplus _{j=1}^{n-1}\st _j}$, which implies that it is linearly independent. 

For any $y\in {\mathcal D}$ fixed, we get that $\{ y\backslash{\overline {\varphi}}(z)\ \mid\  z\in {\displaystyle \bigcup _{j=1}^{n-1}{\mbox {\bf ST}_j}}\}$ is linearly independent in $\st$. 

But, for any $x\in {\mbox {\bf Irr}}$ there exist unique elements $y\in {\mathcal D}$ and $z\in {\mbox {\bf ST}}$ such that $x = y\backslash z$, so we may conclude that 
$$\{{\overline{\eta}}(x)\ \mid x\in {\mbox {\bf Irr}}\} = \{ y\backslash {\overline{\varphi}}(z)\ \mid y\in {\mathcal D}\ {\rm and}\ z\in {\mbox {\bf ST}}\}$$
is linearly independent, which ends the proof.

\bigskip

\bigskip

\section{Tridendriform algebra structures on $\st$}
\medskip

We denote $\st _{qT}$ the $q$-tridendriform algebra, whose underlying vector space is $\K[{\mbox {\bf ST}}]$, which is described in this section. 

The $0$-tridendriform algebra $\st_{0T}$ was introduced by F. Chapoton in \cite{Chapoton}, while the $1$-tridendriform structure $\st_{1T}$ was described in \cite{Loday-Ronco}.
\medskip

For $f\in \ST$, we denote by $s(f)$ the integer $n{-}r$. Using the conventions of Notation \ref{notation3}, we define a $q$-tridendriform structure on $\K[{\mbox {\bf ST}}]$ in terms of stuffles. 
\medskip

\begin{definition} \label{qtridend} {\rm The binary operations $\succ _q$, $\cdot _q$ and $\prec _q$ are defined on the vector space $\K[{\mbox {\bf ST}}]$ as follows:
\begin{enumerate}
\item $x\succ _q y :=\sum _{f\in SH^{\succ}(r,s)} q^{s(f)} f\circ (x\t y)$,
\item $x\cdot _q y :=\sum _{f\in SH^{\bullet}(r,s)} q^{s(f)-1} f\circ (x\t y)$,
\item $x\prec _q y :=\sum _{f\in SH^{\prec}(r,s)} q^{s(f)} f\circ (x\t y)$,\end{enumerate}
for $x\in \ST$ and $y\in {\mbox {\bf ST}_m^s}$, where when $q = 0$ we establish that $q^0 :=1$.}\end{definition}
\medskip

For example, if $f = (2,1,1)\in {\mbox {\bf ST}_3}$ and $g =(1,2)\in {\mbox {\bf ST}_2}$, then
\begin{align}
f \succ_q g &= (2,1,1,3,4) + q (2,1,1,1,3) + q (2,1,1,2,3) + (3,1,1,2,4) + (3,2,2,1,4),\notag\\
f \cdot g &= q (2,1,1,1,2) + (3,2,2,1,3) + (3,1,1,2,3),\notag \\
f \prec g &= q (3,1,1,1,2) + q (3,2,2,1,2) + (4,1,1,2,3) + (4,2,2,1,3) + (4,3,3,1,2).\notag\end{align}
\medskip

\begin{remark} \label{case0} When $q= 0$, the definition of $\succ _0$, $\prec _0$ and $\cdot _0$ are simpler than the general case. For instance, we get that 
$$x\succ _0 y :=\sum _{f\in Sh^{\succ}(r,s)}  f\circ (x\t y),$$
and a similar formula for $\prec _0$. In the case of $\cdot _0$, the sum described in Definition \ref{qtridend} is taken over all $f\in {\mbox {Sh}^{\bullet }(r,s)}$.\end{remark}
\medskip

The following result was proved in \cite{Burgunder-Ronco}, we refer to it for the details of the proof.

\begin{proposition} For any $q\in \K$, the space $\K[{\mbox {\bf ST}}]$ with the operations $\succ _q$, $\cdot _q$ and $\prec_q$ is a $q$-tridendriform algebra.\end{proposition}
\medskip

\begin{notation} \label{notationqtri} {\rm The $q$-tridendriform algebra $(\K[{\mbox {\bf ST}}], \succ _q, \cdot _q, \prec _q)$ described in Definition \ref{qtridend} is denoted  $\st _{qT}$.}\end{notation}
\medskip

The following result is also proved in \cite{Burgunder-Ronco}.

\begin{proposition} The $q$-tridendriform algebra $\st _{qT}$ with the coproduct  $\Delta $, described in Definition \ref{coproduct}, is a conilpotent $q$-tridendriform bialgebra.\end{proposition}

Note that the underlying coalgebra structure of $\st _{qT}$ is $\st$, so the subspace of primitive elements of $\st _{qT}$ is ${\mbox {Prim}(\st )}$.
\bigskip

\bigskip

\section{ $\K[{\mbox {\bf Irr}}]$ as a free ${\mbox {GV}_q}$ algebra}
\medskip

Our final goal is to exhibit a basis of $\st _{qT}$ as a free $q$-tridendriform algebra. The outline of the proof is similar to the one we used to construct a basis of $\st _D$ as a free dendriform algebra. 

In the present section we show that the graded vector space $\K[{\mbox {\bf Irr}}]$ admits a structure of free ${\mbox {GV}_q}$ algebra and in the last section we prove that there exists an isomorphism of ${\mbox {GV}_q}$ algebras from $\K[{\mbox {\bf Irr}}]$ to ${\mbox {Prim}(\st )}$.
\bigskip

\begin{definition} \label{mid} {\rm Let $x = \prod _{j_1<\dots < j_{\lambda (x)}}x\rq \in \ST$ and $y = \prod _{k_1<\dots <k_{\lambda (y)}}y\rq $ be two surjective maps. Define the product
$$x\cdot y =\prod _{j_1<\dots < j_{\lambda (x)} < k_1+r <\dots <k_{\lambda (y)+r}} x\rq \times y\rq.$$}\end{definition}

It is easy to verify that $\cdot $ is associative.
\bigskip

We begin by constructing our basis.

Let $\C_n$ be the set of irreducible elements of ${\mbox {\bf ST}_n}$ defined recursively as follows:\begin{enumerate}
\item $\C_1=\emptyset$.
\item An element $x\in\Irr$ belongs to $\C_n$ if it fulfills one of the following conditions:
\begin{enumerate}[i.]
\item there exist $y\in{\mbox {\bf Irr}_m^s}\setminus \C_{m}$ and $z\in{\mbox {\bf ST}^{n-m}_{r-s}}$ such that $x = y\backslash z$,
\item there exists $x^1\in{\mbox {\bf ST}_m^s}$ and $x^2\in{\mbox {\bf ST}_{n-m}^{r-s+1}}$,  such that  $x = x^1\cdot x^2$.
\end{enumerate}  
\end{enumerate}
\medskip

Note that, if $x= (r,x(2),\dots ,x(n))\in\ST$, then $x = (1)\cdot (x\vert _{\{2,\dots ,n\}})\in {\mathcal C}_n$. And if $x= (x(1),\dots ,x(n-1),r)\in\ST$, then $x = x\vert _{\{1,\dots ,n-1\}}\cdot (1)\in \C_n$.
\medskip

\begin{example}\begin{enumerate}
\item The element $x = (2,5,1,3,5,2,4,5,4)$ belongs to $\C_9$ because it fulfills the third condition. Indeed, $x^1= (2,4,1,3,2)\in{\mbox {\bf ST}_5^4}$, $x^2=(1,2,1)\in{\mbox {\bf ST}_3^2}$, and $x = x^1\cdot x^2$. 
\item The element $y = (4,5,2,3,1)$ belongs $\C_5$ as it verifies the second condition, for $x^1= (3,4,1,2)$ and $x^2=(1)$.  Note that $x1\notin \C_4$.\end{enumerate}
\end{example}
\medskip

\begin{definition} \label{decomposable} An {\it indecomposable map} is an element $x\in \ST$ such that there do not exist surjective maps $x^1$ and $x^2$ satisfying that $x =  x^1\cdot x^2$. We denote by ${\mbox {\bf Indec}}$ the set of indecomposable elements of ${\mbox {\bf ST}}$.\end{definition}

A standard argument shows that for any $x\in {\mbox {\bf ST}}$ there exist a unique integer $p\geq 1$ and unique indecomposable elements $x^1,\dots ,x^p$ such that $x^2,\dots ,x^p$ are irreducible and $x = x^1\cdot \ldots \cdot x^p$. 

For instance, the element $x = (2,3,7,1,3,4,7, 5, 6,7,5)$ may be written as $x = (2,3,5,1,3,4)\cdot (1,2,3,1)$, as $(2,3,7,1,3)\cdot (1,4, 2, 3,4,2)$ or as $(2,3,7,1,3)\cdot (1,2)\cdot (1, 2,3,1)$, but the unique decomposition with $x^2$ irreducible is the first one.
\bigskip

Let $\B_n$ be the set ${\mbox {\bf Irr}_n}\setminus \C_n$, for $n\geq 1$.
\bigskip

\begin{definition} \label{GVonIrr} {\rm We have already introduced the product $\cdot $ on $\K[{\mbox {\bf Irr}}]$.  Define the operations $M_{1n}$ as follows:\begin{enumerate}
\item $M_{1n}(x ; y_1,\dots ,y_n):= x\backslash (y_1\t \dots \t y_n)$, for
any $x \in \bigcup _{n\geq 1}\B_n$ and any family of irreducible functions $y _1,\dots ,  y_n$.
\item for an $x \in \C _n$, we have two possibilities:\begin{enumerate}
\item If $x$ is indecomposable, then $x\in \B _n$ or $x = y\backslash z$, with $y\in \B_m$ and $z\in {\mbox {\bf ST}_{n-m}}$. 

For $x\notin \B_n$, there exist unique irreducible functions $z^1,\dots , z^p$ such that $z = z^1\t \dots \t z^p,$
and $x =M_{1p}(y; z^1,\dots ,z^p)$. 

So, the element $$M_{1n}(x ;w^1,\dots ,w^n)=M_{1n}(M_{1p}(y ; z^1,\dots ,z^p); w^1,\dots , w^n)$$ is well defined applying Definition \ref{def:Gerstenhaber-Voronov} and Lemma \ref{GVbasis}.
\item If $x = x^1\cdot \ldots \cdot x^p$, for some  $x^i\in {\mbox {\bf ST} _{n_i}}$ and some  $p\geq 2$, then, for each $1\leq i\leq p$, we may suppose that either $x ^i\in B_{n_i}$ or $x^i = y^i\backslash z^i$, with $y^i\in \B_{s_i}$.
Again, applying the definition of ${\mbox {GV}_q}$ algebra, $M_{1n}(x ;w_1,\dots ,w_n)$ may be computed in terms of products of elements of type $M_{1l_j}(x^j; w_{i_1},\dots ,w_{i_j})$ and $w_j$.\end{enumerate}\end{enumerate}}\end{definition}
\medskip

So, we have a natural structure of ${\mbox {GV}_q}$ algebra on $\K[{\mbox {\bf Irr}}]$.

\begin{proposition} The algebra $\K[{\mbox {\bf Irr}}]$, with the structure described in Definition \ref{GVonIrr}, is the free ${\mbox {GV}_q}$ algebra spanned by $\bigcup _{n\geq 1}\B_n$.\end{proposition}
\medskip

\begin{proo}
Any  element $x\in{\mbox {\bf Irr}}$ is written uniquely as a product $ x = x^1\cdot \ldots \cdot x^p$, for $p\geq 1$, with $x^i\in \bigcup _{n\geq 1}\B_n$ or $x^i = M_{1n_i}(y^i; z_1^i,\dots , z_{n_i}^i)$, for all $1\leq i\leq k$. So, $\bigcup _{n\geq 1}\B_n$ spans $\K[{\mbox {\bf Irr}}]$ as a free ${\mbox {GV}_q}$ algebra.\end{proo}
 \bigskip
 
 \bigskip
 
 \section{Freeness of $\st _{qT}$ as a $q$-tridendriform algebra}\label{freetridend}
 \medskip
 
 Using   Theorem \ref{idempotents} and arguments similar to those in section \ref{dendbialg}, we shall prove that ${\mbox {Prim}(\st)}$ is the free ${\mbox {GV}_q}$ algebra spanned by $\B$, which implies that $\st _{qT}$ is the free $q$-tridendriform algebra spanned by $\B$.
\medskip 
 
Again, define $\psi^q : \K[{\mbox {\bf Irr}}]\longrightarrow {\mbox {Prim}(\st)}$ as the homomorphism of ${\mbox {GV}_q}$ algebras defined by setting:
$$\psi^q (x) := E(x) = \sum _{\underline l}\alpha _{\underline l} x^{\underline l},$$
for $x\in \B$.
\bigskip

\begin{theorem} \label{principthe} The homomorphism $\psi^q : \K[{\mbox {\bf Irr}}]\longrightarrow {\mbox {Prim}(\st_{qT})}$ is an isomorphism, for all $q\in \K$.\end{theorem}
\bigskip

The rest of this section is devoted to the proof of Theorem \ref{principthe}, which implies that ${\mbox {Prim}(\st)}$ is a free ${\mbox {GV}_q}$ algebra, and therefore that 
$\st _{qT}$ is a free $q$-tridendriform algebra.
\medskip

As $\vert {\mbox {\bf Irr}_n}\vert = {\mbox {dim}_{\K}({\mbox {Prim}(\st)_n}})$, for $n\geq 1$, it suffices to see that $\psi^q \vert _{\K[Irr]_n}$ is surjective, or injective, for all $n\geq 1$.
\bigskip

Consider the following subsets of the set ${\mbox {\bf {Irr}}_{ST}}$:\begin{itemize}
\item ${\mathcal B}$ is the basis.
\item ${\mbox {\bf Br}}$, is the set of elements of the form $x= y\backslash z$, with $y\in {\mathcal B}$, and $z\in {\mbox{\bf ST}}$. The set ${\mbox {\bf Br}}$ is the disjoint union
$${\mbox {\bf Br}}=\bigcup _{n\geq 0}{\mbox {\bf Br}(n)},$$
where ${\mbox {\bf Br}(n)}$ the subset of elements such that $\vert z\vert =n$. For $n=0$, we have ${\mbox {\bf Br}(0)}= {\mathcal B}$.
\item ${\mbox {\bf Prod}}$, is the set of decomposable elements of the form $x=x^1\cdot x^2$, with $x^1$ and $x^2$ irreducible. 
Define $${\mbox {\bf Prod}(p)}:=\{ x\in {\mbox {\bf Prod}}\vert \ {\rm such\ that}\ x = x^1\cdot \ldots \cdot x^p,\ {\rm with}\ x^i\in \bigcup _{l\geq 0}{\mbox {\bf Br}(l)}\},$$ so that
 ${\mbox {\bf Prod}}$ is the disjoint union of $\{ {\mbox {\bf Prod}(p)}\} _{p\geq 2}$. 
\end{itemize} 

Note that, by the definition of $\B$, we get that:\begin{enumerate}
\item ${\mbox {\bf Irr}}$ is the union of ${\mbox {\bf Br}}$ and ${\mbox {\bf Prod}}$,
\item $\B \cap {\displaystyle \bigcup_{m\geq 1}{\mbox {\bf Br}(m)}} = \emptyset$ and $\B\cap {\mbox {\bf Prod}}= \emptyset$.\end{enumerate}

Any element ${\displaystyle x= \prod_{j_1<\dots <j_{\lambda (x)}}x\rq }$ in ${\displaystyle \bigcup_{m\geq 1}{\mbox {\bf Br}(m)}_n}$ satisfies that $x^{-1}(\{1\})\subseteq \{j_{\lambda (x)}+1,\dots ,n\}$.  

On the other hand, if ${\displaystyle x= \prod_{j_1<\dots <j_{\lambda (x)}}x\rq\in {\mbox {\bf Prod}}}$ we have that $x^{-1}(\{1\})\subseteq  \{1,\dots ,j_{\lambda (x)}-1\}$, which implies that 
${\mbox {\bf Br}}\cap {\mbox {\bf Prod}} = \emptyset $.
\bigskip

Let us recall the definition of $\psi ^q(x)$ for $x\in {\displaystyle \bigcup _{l\geq 1} {\mbox {\bf Br}(l)}}$ and $x\in {\mbox {\bf Prod}}$. We denote by $\cdot ^{ST(q)}$ and $M_{1n}^{ST(q)}$ the associative product and the brace operations defined on $\st_{qT}$.
\medskip

\begin{notation}\label{notation1}  As in section \ref{dendbialg}, for $x\in {\mbox {\bf Irr}}$, we denote  $\psi^q (x) = \sum_{a_i\neq 0} a_ix_{i}$.\end{notation}

For $x\in {\mathcal B}_n^r$ and ${\underline l} = (l_1,\dots ,l_p)$, the element $x^{\underline l}$ belongs to ${\mbox {\bf ST}_n^r}$, which implies that 
$\psi^q (x) \in \K[\ST ]$.
\medskip

For example, \begin{itemize}\item $\psi ((2,3,1))= (2,3,1)-(2,1,3)$,\  $\psi ((1,2,1))= (1,2,1) - (1,1,2)$,
\item $\psi ((2,3,4,1))=(2,3,4,1) - (2,3,1,4)$, 
\item $\psi ((2,4,3,1))= (2,4,3,1)- (2,1,4,3) + (2,1,3,4) - (2,3,1,4)$.\end{itemize}
\bigskip

Using Notation \ref{notation3}, we may describe easily $\psi (x)$, for $x\in {\displaystyle \bigcup _{l\geq 1} {\mbox {\bf Br}(l)}}$.

\begin{notation}\label{notation2} {\rm Let $x^1,\dots, x^p$ be a collection of irreducible elements in ${\mbox{\bf ST}}$, with $x^j\in {\mbox{\bf ST}_{n_j}^{r_j}}$. For $q\in K$, we denote by:\begin{enumerate}
\item $\gamma ^{\succcurlyeq _q}(x^1,\dots ,x^p)$ the element:
$$(\dots ((x^1 {\succcurlyeq _q} x^2)\succcurlyeq _q\dots )\succcurlyeq_q x^p=\sum _{f\in {SH}^{\succcurlyeq }(r_1,\dots ,r_p)} q ^{s(f)} f\circ (x^1\times \ldots \times x^p),$$
\item $\gamma ^{\succ_q}(x^1,\dots ,x^p)$ the element:
$$(\dots ((x^1 {\succ_q}x^2)\succ_q\dots )\succ_q x^p=\sum _{f\in {SH}^{\succ}(r_1,\dots ,r_p)} q ^{s(f)} f\circ (x^1\times \dots \times x^p) ,$$
\item $\gamma ^{\prec _q}(x^1,\dots ,x^p)$ the element:
$$x^1\prec_q (x^2 \prec_q  (\dots (x^{p-1}\prec_q x^p)\dots )))= \sum _{f\in {SH}^{\prec}(r_1,\dots ,r_p)}q ^{s(f)} f\circ (x^1\times \ldots \times x^p).$$\end{enumerate}
}\end{notation}

Applying the formula of $M_{1n}^{ST(q)}$ to an element $x\backslash y\in {\mbox {\bf Br}}$, for $x\in {\mathcal B}_n^r$ and $y=y^1\times \ldots \times y^p$, with $y^j\in {\mbox {\bf {Irr}}}$, we get that:
$$\psi ^q(x\backslash y)=\sum _{j=0}^p \bigl ( \sum c_{i,i_1,\dots ,i_p}\gamma ^{\prec _q}(y_{i_1}^1,\dots ,y_{i_j}^j)\succcurlyeq_q\ x^{\underline l}\prec_q \gamma ^{\succcurlyeq _q}(y_{i_{j+1}}^{j+1},\dots ,y_{i_p}^p)\bigr ),\leqno(8.1)$$
where  $\psi^q (x)=\sum_{\underline l} a_{\underline l} x^{\underline l}$, $\psi^q(y^j)=\sum_{b_{i_j}\neq 0} b_{i_j} y_{i_j}^j$, and $c_{i,i_1,\dots ,i_p} \in {\mathbb Z}.$ 
\medskip

Using Notation \ref{notation2}, we rephrase  formula $(8.1)$ as follows:
$$\psi^q (x\backslash y)=\sum _{j=0}^p \bigl ( \sum c_{l,i_1,\dots ,i_p} h\circ (f\circ (y_{i_1}^1\times \ldots \times y_{i_j}^j)\times x^{\underline l}\times g\circ (y_{i_{j+1}}^{j+1}\times \ldots \times y_{i_p}^p))\bigr ),\leqno(8.2)$$
where the sum is taken over all $f\in {\mbox {SH}^{\prec}(r_{i_1},\dots ,r_{i_j})^q}$, $g\in {\mbox {SH}^{\succcurlyeq}(r_{i_{j+1}},\dots ,r_{i_p})^s}$ and $h\in {\mbox {SH}(u , r, s)}$ such that $h(u)\leq h(u+r)$ and $h(u+r) > h(u+r+s)$, for $y_{i_k}^k\in {\mbox {\bf ST}_{n_{i_k}}^{r_{i_k}}}$. 
\bigskip

In a similar way, if $x_1,\dots ,x_p$ is a family of elements in ${\mbox {\bf Br}}$, with $p\geq 2$, and $x = x^1\cdot \ldots \cdot x^p\in {\mbox {\bf Prod}(p)}$, then 
$$\psi^q (x) = \sum \bigl (\sum _{f\in {SH}^{\bullet}(r_{i_1},\dots ,r_{i_p})}c_{i_1\dots i_p} f\circ (x_{i_1}^1\t \dots \t x_{i_p}^p)\bigr ),\leqno(8.3)$$
where $\psi ^q(x^j) = \sum _{a_{i_j}\neq 0}a_{i_j}x_{i_j}^j$, with $x_{i_j}^j\in {\mbox {\bf ST}_{n_j}^{r_{i_j}}}$, for $1\leq j\leq p$.
\bigskip

\bigskip

\centerline{\bf Reduction to the case $q=0$}
\medskip

We want to see that if $\psi ^0$ is an isomorphism, then $\psi^q$ is an isomorphism too, for all $q\in \K$. This result implies that the freeness of $\st _{qT}$ is equivalent to the freeness of $\st _{0T}$, for any $q\in \K$.
\medskip

The proof of the following Lemma is immediate.

\begin{lemma}\label{lemma5} Suppose that $x^1,\dots ,x^p$ is a family of surjective maps such that $x^i\in {\mbox{\bf ST}_{n_i}^{r_i}}$. 
For $f\in {\mbox {SH}(r_1,\dots ,r_p)}$, we get that the map $f\circ (x^1\times \ldots \times x^p)$ belongs to ${\mbox{\bf ST}_{n}^{s}}$, for $s\leq r_1+\dots +r_p$. Moreover, $s= r_1+\dots +r_p$ if, and only if, $f\in {\mbox {Sh}(r_1,\dots ,r_p)}$. In particular, if $f\in {\mbox {SH}^{\bullet}(r_1,\dots ,r_p)}$, then $s < r_1+\dots +r_p$.\end{lemma}
\medskip

For all $x\in {\mbox {\bf Irr}_n^r}$ such that $\psi^q(x) = \sum _{a_i\neq 0}a_ix_i$, the elements $x_i\in {\mbox {\bf ST}_n^{r_i}}$ satisfy that $r_i \leq r$. So, we may  restrict ourselves to work with ${\displaystyle \sum _{{r_i = r}\atop {a_i\neq 0}}a_i x_i}$ instead of $\psi^q(x)$. 
\medskip

Note that for $x\in \B _n^r$, the element $\psi ^q(x) = \sum _{\underline l} \alpha_{\underline l} x^{\underline l}$, with $x^{\underline l}\in \ST$, does not depend on $q$. 
\medskip

Suppose that $x = y\backslash (z^1\t\dots \t z^p)$, with $y\in \B _m^s$, $z^j\in {\mbox {\bf Irr}_{n_j}^{l_j}}$ such that $\psi^q(z^j) = \sum _{b_{i_j}\neq 0}b_{i_j}z_{i_j}^j$, for $1\leq j\leq p$. By a recursive argument, we may assume that:
$$\psi ^0(z_j) = \sum _{{l_{i_j} = l_j}\atop{b_{i_j}\neq 0}}b_{i_j}z_{i_j}^j,$$
for $1\leq j\leq p$. So, ${\displaystyle \sum _{{r_i = r}\atop {a_i\neq 0}}a_i x_i }$ is taken over all $x_i$ satisfying that
$$x_i = h\circ (f\circ (z_{i_1}^1\times \ldots \times z_{i_k}^k)\times y^{\underline l}\times g\circ (z_{i_{k+1}}^{k+1}\times \ldots \times z_{i_p}^p)),$$
with $z_{i_j}^j\in {\mbox {\bf ST}_{n_j}^{l_j}}$, for $1\leq j\leq p$, 
$f\in {\mbox {Sh}^{\prec}(l_1,\dots ,l_k)}$, $g\in {\mbox {Sh}^{\succ}(l_{k+1},\dots ,l_p)}$, and $h\in {\mbox {Sh}(l_1+\dots +l_k,s,l_{k+1}+\dots +l_p)}$ such that 
$h(l_1+\dots +l_k) < h(l_1+\dots +l_k +s)$ and $h(l_1+\dots +l_k +s) > h(l_1+\dots +l_p +s)$. The coefficient of $x_i$ in $\psi ^q(x)$ being $\alpha_{\underline l}b_{i_1}^1\cdot \dots \cdot b_{i_p}^p$. We get that ${\displaystyle \sum _{{r_i = r}\atop {a_i\neq 0}}a_i x_i =\psi^0(x)}$.
\medskip

In a similar way, if $x = x^1\cdot \ldots \cdot x^p$, with $x^j\in {\mbox {\bf Br}_{n_j}^{r_j}}$, $1\leq j\leq p$, then the unique elements $x_{i}$ which belong to $\ST$  are of the form
$$x_{i} = f\circ (x_{j_1}^1\t \dots \t x_{j_p}^p),$$
 with $x_{j_k}^k\in  {\mbox {\bf ST}_{n_k}^{r_k}}$, for $1\leq k\leq p$, and $f \in {\mbox {Sh}^{\bullet }(r_1,\dots ,r_p)}$. 
 
As $n_j < n$, we assume that $\psi^0(x^j) = {\displaystyle \sum _{{r_{i_jj}=r_j}\atop {b_{i_j}\neq 0}}b_{i_j}x_{i_j}^j}$, for $1\leq j\leq p$.
So, $$ \sum _{{r_i=r}\atop {a_i\neq 0}}a_ix_i=
\sum b_{i_1}\cdot \dots \cdot b_{i_p} \bigl (\sum _{f}f\circ (x_{j_1}^1\t \dots \t x_{j_p}^p)\bigr),$$
where $\psi ^0(x^j) =\sum _{b_{i_j}\neq 0} b_{i_j} x_{i_j}^j$ and $f \in {\mbox {Sh}^{\bullet }(r_1,\dots ,r_p)}$, 
which implies that ${\displaystyle \sum _{{r_i=r}\atop {a_i\neq 0}}a_ix_i= \psi ^0(x)}$.

\medskip

The following Proposition is an immediate consequence of the above arguments and of Remark \ref{case0}.

\begin{proposition} \label{proposition1}  If $\{\psi^0(x)\ \mid\ x\in {\mbox {\bf Irr}_n}\}$ is linearly independent in ${\mbox{Prim}(\st)}$, then $\{ \psi^q (x)\ \mid\ x\in {\mbox {\bf Irr}_n}\}$ is linearly independent in ${\mbox{Prim}(\st)}$ too, for all $n\geq 1$ and all $q\in \K$.\end{proposition}
\medskip

Note that, even if we have not specified it, the products $\cdot ^{ST}$ and $M_{1n}^{ST}$ depend on the tridendriform structure $\st _{qT}$ considered. From now on, as we have proved that it suffices to work with $q=0$, the operations $\cdot ^{ST}$ and $M_{1n}^{ST}$ will be the ones defined in $\st _{0T}$.
\bigskip

\bigskip

\centerline{\bf Reduction ${\mathbb M}$}
\medskip

In Lemma \ref{lemmaMalphaE} we proved that for any $x\in \ST$, and any family ${\underline l} = (l_1,\dots ,l_p)$, with $0 = l_0 < l_1<\dots <l_p < r$, we have that ${\mathbb M}(x^{\underline l}) \leq {\mathbb M}(x)$ for the lexicographic order. 

We want to prove that, for any $x\in {\mbox {\bf Irr}}$ such that $\psi^0(x) = \sum _{a_i\neq 0}a_ix_{i}$, we have that 
${\mathbb M}(x_{i})\leq {\mathbb M}(x)$ for the lexicographic order.
\medskip

Clearly, the result holds for $x\in \B$, since $\psi ^0(x) =  \sum _{\underline l}\alpha _{\underline l} x^{\underline l}$.
\medskip

The proof of the following Lemma is easily obtained by mimicking the proof of Lemma \ref{lemma4}, because the arguments we used to prove it still apply when we replace $\eta$ by $\psi ^0$.

\begin{lemma} \label{lemma6} Let $x\in {\mbox {\bf Irr}}$ such that $\psi ^0(x) = \sum _{a_i\neq 0} a_i x_{i}$. Any element $x_{i}$ satisfies that ${\mathbb M}(x_{i})\leq {\mathbb M}(x)$ for the lexicographic order. Moreover, for $x = y\backslash (z^1\t \dots \t z^p)$, with $y\in \B_m^s$ and $z^1,\dots ,z^p$ irreducibles, if ${\mathbb M}(x_{i}) = {\mathbb M}(x)$, then 
$$ x_{i} = h\circ (y^{\underline l}\t (g\circ (z_{i_1}^1\t \dots \t z_{i_p}^p)),$$
for some family ${\underline l}$ such that ${\mathbb M}(y^{\underline l}) = {\mathbb M}(y)$, $g\in {\mbox {Sh}^{\succ}(r_1,\dots ,r_p)}$ and $h\in {\mbox {Sh}^{\leq}(s,r_1+\dots +r_p)}$, 
where $z^k\in {\mbox {\bf Irr}_{n_k}^{r_k}}$
is such that $\psi ^0(z^k) =\sum _{b_{i_k}\neq 0} b_{i_k}z_{i_k}^k$, for $1\leq k\leq p$.\end{lemma}
\bigskip

For an element $x\in {\mbox {\bf Prod}}$, we have a similar result.

\begin{lemma} \label{lemma7} Let $x = x^1\cdot \ldots \cdot x^p$ such that $x^i\in {\mbox {\bf Br}_{n_i}^{r_i}}$, for $1\leq i\leq p$. If $\psi ^0(x^j) =\sum _{a_{i_j}\neq 0} a_{i_j} x_{i_j}^j$ and $f \in {\mbox {Sh}^{\bullet}(r_1,\dots ,r_p)}$, then $${\mathbb M}(f\circ (x_{i_1}^1\t \dots \t x_{i_p}^p)) \leq {\mathbb M}(x).$$ The equality holds if, and only if ${\mathbb M}(x_{i_j}^j) = {\mathbb M}(x^j)$, for $1\leq j\leq p$.\end{lemma}\medskip

\begin{proo} Denote $w = f\circ (x_{i_1}^1\t \dots \t x_{i_p}^p)$. It is easily seen that
$${\mathbb M}(w)_j = \begin{cases} {\mathbb M}(x_{i_p}^p), &\ {\rm for}\ j = \lambda (w),\\
{\mathbb M}(x_{i_k}^k)_l, &\ {\rm for}\ j = \lambda(x_{i_1}^1)+\dots +\lambda(x_{i_{k-1}}^{k-1}) + l,\\
{\mathbb M}(x_{i_k}^k)_{\lambda (x_{i_k}^k)} + {\mathbb M}(x_{i_{k+1}}^{k+1}), &\ {\rm for}\ j=\lambda(x_{i_1}^1)+\dots +\lambda(x_{i_k}^k).\end{cases}$$

The result follows, by applying that ${\mathbb M}(x_{i_k}^k) \leq {\mathbb M}(x^k)$, $1\leq k\leq p$, and using Lemmas \ref{lemmaMalphaE} and \ref{lemma6}.\end{proo}

So, we get a second reduction. Let 
$$\psi_{\mathbb M}(x) := \sum _{{\mathbb M}(x_{i}) = {\mathbb M}(x)}a_i x_{i},$$
where $\psi ^0(x) = \sum _{a_i\neq 0}a_ix_{i}$.

\begin{proposition}\label{proposition2} If the set $\{ \psi _{\mathbb M}(x)\mid x\in {\mbox {\bf Irr}}\}$ is linearly independent in $\st $, then $\{ \psi ^0(x)\mid x\in {\mbox {\bf Irr}}\}$ is linearly independent in ${\mbox {Prim}(\st)}$.\end{proposition}
\bigskip

\bigskip

\centerline{\bf The image of $\Br$ under $\psi_{\mathbb M}$}
\medskip

We want to prove that $\psi ^0(\Br(l))\subseteq \K[{\mbox {\bf Red}}\cup {\displaystyle \bigcup _{i=0}^l\Br(i)}]$, where ${\mbox {\bf Red}}$ denotes the set of reducible elements. Our results are similar to the ones of section \ref{dendbialg}, but in the tridendriform case.
\medskip

\begin{lemma} \label{lemma8} Let $x\in \B_n^r$, $y\in {\mbox {\bf ST}_m^s}$ and $f\in {\mbox {Sh}(r,s)}$. If $f \neq \epsilon(r,s)$, then 
$$f \circ (x\times y)\in \bigcup_{l=0}^{m-1}\Br(l).$$
\end{lemma}
\medskip

\begin{proo} Suppose that $x = \prod_{j_1<\dots <j_{\lambda(x)}}x\rq$ and $f = \prod_{r}f\rq$, with $f\rq \in {\mbox {Sh}(r{-}1,s)}$. We have that
$f \circ (x\times y) = \prod_{j_1<\dots <j_{\lambda(x)}}f\rq \circ (x\rq \t y)$.

If $f \circ (x\times y)$ is reducible, then $f \circ (x\times y) = h^1\t h^2$, with $\vert h^1\vert < j_1$. So,
$$x = f \circ (x\times y)\vert_{\{1,\dots ,n\}} = h^1\t  (h^2\vert _{\{1,\dots ,n-\vert h^1\vert\}}),$$
which is false since $x$ is irreducible. 
\medskip

If $f \circ (x\times y)\in \Prod$, then $f\rq \circ (x\rq \t y)= h^1\times h^2$, with $\vert h^1\vert < j_{\lambda (x)}{-}\lambda(x)$.
So, $$x = f \circ (x\times y)\vert_{\{1,\dots ,n\}} = \prod_{j_1<\dots <j_{\lambda(x)}}h^1\t  (h^2\vert _{\{1,\dots ,n{-}\vert h^1\vert\}}),$$
which does not happen because $x\notin \Prod$. 
\medskip

Suppose that $f \circ (x\times y)\in \Br(l)$, that is $f \circ (x\times y) = z\backslash w$, with $z\in \B$ and $w\in {\mbox {\bf ST}}$ with $\vert w\vert = l$. Again, we have that 
$$x =  f \circ (x\times y)\vert_{\{1,\dots ,n\}} = (z\backslash w)\vert _{\{1,\dots ,n\}},$$
which implies that $\vert z\vert \geq n$, otherwise, $x = z\backslash (w\vert _{\{1,\dots ,n{-}\vert z\vert\}})$ is not in $\B$.

We want to see that $\vert z\vert > n$. If $\vert z\vert = n$, then 
$$f \circ (x\times y)(i) = \begin{cases} f(x(i)) = z(i) + p,& {\rm for}\ 1\leq i\leq n\\
f(y(i {-} n)+r) = w(i {-} n),& {\rm for}\ n +1\leq i\leq n + m,\end{cases}$$ 
where $w\in {\mbox {\bf ST}_l^p}$. 

But $x (\{1,\dots ,n\}) = \{1,\dots ,r\}$ and $y(\{1,\dots ,m\}) + r = \{ r+1,\dots ,r+s)\}$, which implies that 
$$f(\{1,\dots ,r\}) = \{p+1,\dots ,p+r\}\ {\rm and}\ f(\{r+1,\dots ,r+s\}) = \{1,\dots ,p\}.$$
 As the unique $f\in {\mbox {Sh}(r,s)}$ satisfying this condition is $f = \epsilon (r,s)$, the proof is over.\end{proo} 
\medskip

\begin{lemma}\label{lemma9} Let $x\in {\mbox  {\bf Red}_n^r}$, $y\in {\mbox {\bf ST}_m^s}$ and $f\in {\mbox {Sh}^{\prec }(r,s)}$. The surjection
$f\circ (x\times y)$ belongs to ${\mbox {\bf Red}} \cup {\displaystyle \bigcup _{l=0}^{m-1}{\mbox {\bf Br}(l)}}$.\end{lemma}
\medskip

\begin{proo} We have to check that $f\circ (x\times y)\notin {\mbox {\bf Prod}} \cup {\displaystyle\bigcup _{l\geq m}{\mbox {\bf Br}(l)}}$. 
\medskip

Let $x=x^1\times x^2$, for $x^1\in {\mbox {\bf ST}_{n_1}^{r_1}}$ and $x^2 \in {\mbox {\bf Irr}_{n_2}^{r_2}}$ such that $x^2=\prod _{j_1<\dots <j_{\lambda(x^2)}}x^{2\rq}$. We have that
$$f\circ (x\times y)=\prod _{j_1<\dots <j_{\lambda (x^2)}} f\rq \circ (x^1\times x^{2\rq} \times y),$$
with $f\rq \in {\mbox {Sh}(r {-} 1,s)}$.
\medskip

Suppose that $f\circ (x\times y)\in \Prod $. In this case, $f\rq \circ (x^1\times x^{2\rq} \times y)=h^1\times h^2$, with $ \vert h^1\vert < j_{\lambda (x^2)} {-} k$, which implies that $f\rq (i)<f\rq (j)$, for $i\leq j_1$ and all $j\geq r$. Using that $f\rq (1)<\dots <f\rq (r{-}1)$, we get that $f \circ (x\times y)$ is reducible. So, $f\circ (x\times y)\notin {\mbox {\bf Prod}}$.
\medskip

If $f\circ (x\times y)= z\backslash w$, with $z\in {\mathcal B}$, then $\vert z\vert >n$ because all the elements of the set $\{ f(x(1)),\dots ,f(x(n_1))\}$ are smaller than all the elements of the set $\{ f(x(n_1+1)),\dots ,f(x(n))\}$, and therefore $\vert w\vert < m$, which ends the proof. \end{proo}
\bigskip

Let $x\in \B_n^r$ and $y=y^1\times \ldots \times y^p\in {\mbox {\bf ST}_m^s}$, with $y^j\in {\mbox {\bf Irr}_{m_j}^{s_j}}$, $1\leq j\leq p$. Suppose that
$\psi^0 (y^j)=\sum_{b_k^j\neq 0} b_k^jy_{k}^j$, with $y_{k_j}^j \in {\mbox {\bf ST}_{m_j}^{s_j}}$.
\medskip

We proved in Lemma \ref{lemma6} that, for $\psi^0 (x\backslash y) =\sum _{c_i\neq 0} c_i w_{i}$, the unique elements $w_{i}$ such that ${\mathbb M}(w_{i}) = {\mathbb M}(x\backslash y)$ are of the form
$$w_{i} = f \circ (x^{\underline l}\times g\circ (y_{k_1}^1\times \ldots \times y_{k_p}^p)),$$ 
with ${\mathbb M}(x^{\underline l}) = {\mathbb M}(x)$, for some $f\in {\mbox {Sh}^{\prec}(r,s)}$ and $g\in {\mbox {Sh}^{\succ}(s_1,\dots ,s_p)}$.

As $x\in \B$, either $x^{\underline l}=x$ or $x^{\underline l}$ is reducible. In the second case, applying Lemma \ref{lemma9}, we get that  $w_{i}\in {\mbox {\bf Red}} \cup {\displaystyle \bigcup _{l=0}^{m-1}{\mbox {\bf Br}(l)}}$. 

On the other hand, if $f\neq \epsilon_{r,s}$, Lemma \ref{lemma8}  implies that the element $w_{i}$ belongs to ${\mbox {\bf Red}} \cup\ {\displaystyle \bigcup _{l=0}^{m-1}{\mbox {\bf Br}(l)}}$. 
\medskip

The argument above proves the following result.

\begin{proposition} \label{proposition3} Let $x$ be an element in ${\mathcal B}_n$ and let $y = y^1\times \ldots \times y^p$ be a surjective map in ${\mbox {\bf ST}_n}$, with $y^1,\dots ,y^p$ irreducibles. Suppose that $\psi ^0 (y^j) = \sum_{b_{k_j}^j\neq 0} b_{k_j}^jy_{k_j}^j$ .

If $\psi _{\mathbb M} (x\backslash y) = \sum _{c_i\neq 0}c_i w_{i}$, then an element $w_{i}$ which is not of the form
$$w_{i} = x\backslash (g\circ (y_{k_1}^1\times \ldots \times y_{k_p}^p)),$$
for some $g\in {\mbox {Sh}^{\succ}(s_1,\dots ,s_p)}$, satisfies that $w_{i} \in {\mbox {\bf Red}} \cup\ {\displaystyle \bigcup _{l=0}^{m-1}{\mbox {\bf Br}(l)}}$.\end{proposition}
 \bigskip

\begin{definition} \label{defnpsibr} Let $x = y\backslash z$ , with $y\in \B_n$ and $z = z^1\t \dots \t z^p\in {\mbox {\bf ST}_m}$, such that $z^j\in {\mbox {\bf Irr}_{m_j}}$. Define
$$\psi _{Br}(x) := y\backslash \omega ^{\succ}(\psi ^0(z^1),\dots ,\psi ^0(z^p)) = \sum c_{i_1\dots i_k}(\sum _{f\in Sh^{\succ }(s_1,\dots ,s_j)}y\backslash f\circ (z_{i_1}^1
\t\dots \t z_{i_p}^p),$$
where $\psi ^0(z^j)= \sum _{a_{i_j}\neq 0}a_{i_j}z_{i_j}^j$, and $c_{i_1\dots i_k} = a_{i_1}\dots a_{i_p}$.\end{definition}
\medskip

Propositions \ref{proposition2} and \ref{proposition3} imply the following result:

\begin{corollary} \label{corollarybrace} \begin{enumerate}\item The image of $\Br$ under $\psi _{\mathbb M}$ is a subspace of $\K[{\mbox {\bf Red}}\cup \Br ]$.
\item If the set $\{ \psi _{Br}(x)\mid x\in \Br\}$ is linearly independent in $\K[{\mbox {\bf Red}}\ \cup \Br ]$, then $\{ \psi ^0(x)\mid x\in \Br\}$ is linearly independent in ${\mbox {Prim}(\st)}$.\end{enumerate}\end{corollary}
 \bigskip
 
 \bigskip

 \centerline{\bf The image of ${\mbox {\bf Prod}}$ under $\psi_{\mathbb M}$}
\bigskip

Given elements $x^i = \prod _{j_1^i<\dots  <j_{\lambda(x^i)}^i}x^{i\rq} \in {\mbox {\bf ST}_{n_i}^{r_i}}$, for $1\leq i\leq p$, recall that 
$$x^1\cdot \ldots \cdot x^p =  {\displaystyle  \prod _{j_1^1<\dots <j_{\lambda (x^p)}^p+n_1+\dots +n_{p-1}} (x^{1\rq}\times \ldots \times x^{p\rq})}.$$

For $p\geq 2$, ${\mbox {\bf Prod}(p)}$ is the set of elements $x = x^1\cdot \ldots \cdot x^p\in {\mbox {\bf Irr}}$ such that $x^i\in {\mbox {\bf Br}}$, $1\leq i\leq p$.

We want to give a description of the image 
$$\psi_{\mathbb M}(x^1\cdot ^{ST}\ldots \cdot ^{ST} x^p)= \psi_{\mathbb M}(x^1)\cdot ^{ST}\ldots \cdot ^{ST} \psi_{\mathbb M}(x^p),$$
where $x^1\cdot ^{ST}\ldots \cdot ^{ST} x^p =\sum _f  {\displaystyle  \prod _{j_1^1<\dots <j_{\lambda (x^p)}^p+n_1+\dots +n_{p-1}} f\circ (x^{1\rq}\times \ldots \times x^{p\rq})}$, and the sum is taken over all  $f \in {\mbox {Sh}(r_1{-}1,\dots ,r_p{-}1)}$.
\bigskip

\bigskip

The following result completes the ones of Lemmas \ref{lemma8} and \ref{lemma9}.

 \begin{lemma}\label{lemma10} Let $x\in {\mbox {\bf Indec}}$ and $y\in {\mbox {\bf ST}_m^s}$ be surjective maps, and let $f \in {\mbox {Sh}^{\prec}(r,s)}$. The element $f\circ (x\t y)$ is indecomposable.\end{lemma}
 \medskip
 
 \begin{proo} The result has been proved for $x\in {\mbox {\bf Red}} \cup {\mathcal B}$. Suppose that $x = \prod _{j_1<\dots <j_{\lambda (x)}}x\rq $ is indecomposable and $x\notin {\mbox {\bf Red}} \cup {\mathcal B}$. We have that:\begin{enumerate}
 \item $1 < j_1$ and $\lambda (x) < n$,
 \item if $x\rq = x^1\times x^2$, then $x^1$ is irreducible and $\vert x^1\vert > j_{\lambda (x)} {-} \lambda (x)$,
 \item $f \circ (x\times y)= \prod _{j_1<\dots <j_{\lambda (x)}}f \rq \circ (x\rq \times y),$
 with $f \rq \in {\mbox {Sh}(r {-}1,s)}$. \end{enumerate}
 
Note that $f\circ (x\times y)\notin  {\mbox {\bf Indec}}$ if, and only if, $f\rq \circ (x\rq \times y)= z^1\times z^2$ with $j_1{-}1<\vert z^1\vert \leq j_{\lambda (x)} {-} \lambda (x)$. 

As $f\rq (1) < \dots < f\rq (r{-}1)$, we get that
$$\displaylines {
x\rq ={\mbox {std}( f\rq \circ (x\rq \times y)\vert _{\{1,\dots ,n{-}\lambda (x)\}})}=\hfill\cr
\hfill {\mbox {std}(z^1\times z^2)\vert _{\{1,\dots ,n{-}\lambda (x)\}}}= z^1\times {\mbox {std}(z^2\vert _{\{1,\dots ,n + m {-}\lambda (x) {-} \vert z^1\vert \}})},\cr }$$
which implies $\vert z_1\vert \geq n+m {-} \lambda (x) > j_{\lambda (x)} {-}\lambda (x)$, in contradiction with the fact that $\vert z^1\vert \leq j_{\lambda (x)} {-} \lambda (x)$. So, $f\circ (x\times y)$ is indecomposable.\end{proo} 
\medskip

\begin{proposition} \label{braceindec} For any surjective map $x\in {\mbox {\bf Br}}$, the element $\psi_{\mathbb M}(x)$ belongs to $\K[{\mbox {\bf Indec}}]$.\end{proposition}
\medskip

\begin{proo} Suppose that $x = \prod _{j_1<\dots <j_{\lambda (x)}} x'\in {\mathcal B}_n^r$,  and that ${\underline l} = (l_1,\dots ,l_p)$ is a family of positive integers $1\leq l_1<\dots <l_p<r$. The element $x^{\underline l}$ is decomposable if, and only if,  $x\vert ^{\{l_p+1,\dots ,r\}}$ is decomposable. 

As $x\in {\mbox {\bf Indec}}$, we get that $x^{-1}(\{1,\dots, l_p\}) \cap \{j_1+1,\dots ,n\}\neq \emptyset $, and therefore Lemma \ref{lemmaMalphaE} implies that ${\mathbb M}(x^{\underline l}) < {\mathbb M}(x)$.
\medskip

Let $x = y\backslash z$ be an element in ${\mbox {\bf Br}(l)}$, for $l\geq 1$, and suppose that $\psi_{\mathbb M}(x) = \sum _{a_i\neq 0}a_ix_i$. For any $i$, we know that
$$x_i = h\circ (y^{\underline l} \times {\overline z}),$$
for some ${\underline l}$ such that ${\mathbb M}(y^{\underline l}) = {\mathbb M}(y)$, some element ${\overline z}\in {\mbox{\bf ST}}$ and some $h\in {\mbox {Sh}^{\prec}(s, r-s)}$, where $y\in {\mathcal B}_m^s$. 

As $y^{\underline l}\in {\mbox {\bf Indec}}$, applying Lemma \ref{lemma10} we get that $x_i$ is indecomposable.\end{proo}
\bigskip

\begin{remark} \label{elementsPL} Let $x = \prod _{j_1<\dots <j_{\lambda (x)}}x\rq$ be an indecomposable surjection, there exist unique elements $y, z$ and $w$ such that:\begin{enumerate}
\item $y, w\in {\mbox {\bf ST}}\cup \{1_{\K}\}$  and $z\in {\mbox {\bf Irr}}$,
\item $x =   \prod _{j_1<\dots <j_{\lambda(x)}} (y\t z\t w)$,
\item $0\leq \vert y\vert < j_1{-}1$ and $0\leq \vert w\vert < n{-} j_{\lambda (x)}$.\end{enumerate}
If $x$ is irreducible, then $y =1_{\K}$.\end{remark}
\medskip

\begin{example} Let $x = (2, 1, 5, 7, 3, 7, 3, 4, 5, 6)$, we have that $y=(2,1)$, $z= (3, 1,1, 2, 3)$ and $w = (1)$.\end{example}
\medskip

\begin{remark} \label{genBruhatorder} Note that if $f = \prod_{r_1<\dots < r_1+\dots +r_p} f\rq $ and $g =\prod_{r_1<\dots < r_1+\dots +r_p} g\rq $ are two elements in ${\mbox {Sh}^{\bullet}(r_1,\dots ,r_p)}$, then $f < g$ for the weak Bruhat order in ${\mbox {\bf ST}}$ if, and only if, $f\rq < g\rq$ for the weak Bruhat order on $S_{r_1+\dots +r_p{-}p}$.\end{remark}

For any positive integer $\lambda$ and any ${\overline m}= (m_1,\dots ,m_p)$, $p\geq 1$, we denote by ${\mbox {\bf ST}_{\lambda, {\overline m}}}$ the set of $x\in {\mbox {\bf ST}}$ such that $\lambda (x) = \lambda $ and ${\mathbb M}(x) = {\overline m}$. 
Define the order $\leq _{wB}$ on ${\mbox {\bf Prod}(p)_{\lambda, {\overline m}}}$ as the transitive relation spanned by:
$$g\circ (x^1\t \ldots \t x^p) \leq _{wB} f\circ (x^1\t \ldots \t x^p),$$
for any pair of permutations $f$ and $g$ in ${\mbox {Sh}^{\bullet}(r_1,\dots ,r_p)}$ such that $f \leq g$ for the weak Bruhat order, where $x^j\in {\mbox {\bf Br}_{n_j}^{r_j}}$, $1\leq j\leq p$.
\medskip

The order $\leq _{wB}$ is well defined on ${\mbox {\bf Prod}(p)_{\lambda, {\overline m}}}$ by Proposition \ref{weakBST}. 
\bigskip

 \begin{proposition} \label{proposition4} Let $x^1, \dots , x^p$ be a collection of maps, with $x^i\in {\mbox {\bf Indec}_{n_i}^{r_i}}$ for $1\leq i\leq p$, and let $f \in {\mbox {Sh}^{\bullet}(r_1,\dots ,r_p)}$. \begin{enumerate}
 \item If there exist surjections $y^1, \ldots, y^q$ such that   $f \circ (x^1\t \ldots \t x^p)= y^1\cdot \ldots \cdot y^q$, then $q\leq p$. 
 \item If $f \circ (x^1\t \ldots \t x^p)= y^1\cdot \ldots \cdot y^p$, with $y^j\in {\mbox {\bf Irr}}$ for $2\leq j\leq p$, then one of the following conditions is satisfied:\begin{enumerate}
 \item ${\displaystyle f = \prod_{r_1<\dots < r_1+\dots +r_p} 1_{r_1+\dots +r_p{-}p}}$ and $x^1\cdot \ldots \cdot x^p = y^1\cdot \ldots \cdot y^p$,
 \item ${\displaystyle \prod_{r_1<\dots < r_1+\dots +r_p} 1_{r_1+\dots +r_p{-}p} }< f$ for the weak Bruhat order, and therefore $x^1\cdot \ldots \cdot x^p < y^1\cdot \ldots \cdot y^p$ for the order $<_{wB}$,
 \item ${\displaystyle f = \prod_{r_1<\dots < r_1+\dots +r_p} 1_{r_1+\dots +r_p{-}p}}$ and $(\vert y^p\vert, \dots ,\vert y^1\vert ) < (\vert x^p\vert , \dots ,\vert x^1\vert )$ for the lexicographic order.\end{enumerate}\end{enumerate}
  \end{proposition}
\medskip

\begin{proo} For ${\displaystyle x^i=\prod _{j_1^i<\dots <j_{\lambda(x^i)}^i}x^{i\rq }}$ and ${\displaystyle f _{r_1<\dots < r_1+\dots +r_p}=\prod  f\rq }$, we get that
$$f \circ (x^1\t \ldots \t x^p)=\prod _{j_1^1<\dots <j_{\lambda(x^1)}^1<j_1^2+n_1< \dots <j_{\lambda(x^p)}^p+n_1+\dots +n_{p-1}} f\rq \circ (x^{1\rq} \t \ldots \t x^{p\rq}).$$
\medskip

We proceed by induction on $p$.
For $p=1$, the result is immediate.
\medskip

Suppose that $p\geq 2$, and that $y^j=\prod _{l_1^j<\dots <l_{\lambda (y^j)}^j}y^{j\rq} \in {\mbox {\bf ST}_{m_j}^{s_j}}$, for $1\leq j\leq q$.  
\bigskip

$I)$ If $m_q  \leq n_p$, then $\lambda (y^q)\leq \lambda (x^p)$. We have to consider two cases: \begin{enumerate}
\item if $\lambda (y^q) < \lambda (x^p)$, there exists $1< k_0\leq n {-} m_q$ such that:
$$x^p= (y^1\cdot \ldots \cdot y^{q{-}1})\vert _{\{k_0,\dots ,n{-}m_q\}}\cdot y^q,$$
which is impossible because $x^p$ is indecomposable.
\item if $\lambda (y^q) = \lambda (x^p)$, then either $x^p = y^q$, or there exists $z\in {\mbox {\bf ST}_{n_p{-}m_q}}$ such that $x^p = z\times y^q$.

Note that if $x^p$ is irreducible, then the unique possibility is $x^p = y^q$, but we only assume that $x^p$ is indecomposable. 
\medskip

If $x^p = y^q$, then $y^q\in {\mbox {\bf ST} _{n_p}^{r_p}}$, 
$$f_1\circ (x^1\times \ldots \times x^{p{-}1}) = y^1\cdot \ldots \cdot y^{q{-}1},$$ 
where $f_1:= f\vert _{\{1,\dots ,r_1+\dots +r_{p{-}1}\}} \in {\mbox {Sh}^{\bullet}(r_1,\dots ,r_{p{-}1})}$. Applying a recursive argument, we get that:\begin{enumerate}
\item $q\leq p$,
\item when $q = p$,\begin{enumerate}\item  if $1_{r_1+\dots +r_{p{-}1}{-}p+1} < f_1\rq $, then 
$x^1\cdot \ldots \cdot x^p < y^1\cdot \ldots \cdot y^p$, because
$1_{r_1+\dots +r_p{-}p} < f\rq $.
\item if $f_1\rq = 1_{r_1+\dots +r_{p{-}1}{-}p+1}$, then 
$(\vert y^{p-1}\vert ,\dots ,\vert y^1\vert ) < (\vert x^{p-1}\vert,\dots ,\vert x^1\vert )$. 
So, we get that $f\rq = 1_{r_1+\dots +r_p{-}p}$ and $(\vert y^p\vert, \dots ,\vert y^1\vert ) < (\vert x^p\vert , \dots ,\vert x^1\vert )$.\end{enumerate}\end{enumerate}
\medskip

Otherwise, suppose $x^p=z\times y^q$ is reducible, with $\vert z\vert >0$. We get that $f\rq = f\rq\vert _{\{1,\dots ,r{-}s_p{-}p+1\}} \t 1_{s_p{-}1}.$
 There exists a unique way to write down:
$$f\rq\vert _{\{1,\dots ,r{-}s_p{-}p+1\}} = f _2\circ (1_{r_1+\dots +r_{p{-}2}{-}p+2}\times f _3),$$
for a pair of surjective maps $f _2\in {\mbox {Sh}(r_1{-}1,\dots ,r_{p{-}2}{-}1,r_{p{-}1}+r_p{-}s_q{-}1)}$ and $f _3\in {\mbox {Sh}(r_{p{-}1}{-}1,r_p{-}s_q)}$. 
\medskip

Define the element ${\tilde {x}^{p{-}1}} := (\prod_{r_{p{-}1}}f _3) \circ (x_{p{-}1}\times z)$. We get that:
$$(\prod_{r_1<\dots <r_1+\dots +r_{p-1}} f _2) \circ (x^1\times \ldots \times x^{p{-}2}\times {\tilde {x}^{p{-}1}})= y^1\cdot \ldots \cdot y^{q{-}1},$$
with ${\displaystyle \prod_{r_1<\dots <r_1+\dots +r_{p-1}} f_2\in {\mbox {Sh}^{\bullet}(r_1,\dots ,r_{p{-}2},r_{p{-}1}+r_p{-}s_p)},}$
and 

\noindent ${\displaystyle \prod_{r_{p{-}1}}f _3\in {\mbox {Sh}^{\prec}(r_{p{-}1}, r_p{-}s_p)}}$. \end{enumerate}

From Lemma \ref{lemma10}, we get that  ${\tilde {x}^{p{-}1}}\in {\mbox {\bf Indec}}$. A recursive argument states that $q\leq p$. 
\medskip

If $q = p$, as $\vert y^p\vert < \vert x^p\vert$,  we get that $(\vert y^p\vert, \dots ,\vert y^1\vert ) < (\vert x^p\vert , \dots ,\vert x^1\vert )$.
\bigskip

$II)$  Suppose that $n_p < m_q$. Let $k_0 \leq p{-}1$ be the minimal integer such that $n_1+\dots + n_{k_0} > n {-} m_q$. We have to consider two cases:\begin{enumerate}
\item when $m_q = n_{k_0}+\dots +n_p$, we get that 
$$f\vert _{\{1,\dots ,r_1+\dots +r_{k_0{-}1}\}}\circ (x^1\t \ldots \t x^{k_0{-}1}) = y^1\cdot \ldots \cdot y^{q{-}1}.$$

\noindent Applying a recursive argument, we get that $q{-}1 \leq k_0{-}1 < p{-}1$, so $q < p$ and the result is proved.
\item when $m_q = l+n_{k_0+1}+\dots +n_p$, for $1\leq l< n_{k_0}$, as $x^{k_0}$ is indecomposable we get that there $x^{k_0} = \prod _{j_1^{k_0}<\dots <j_{\lambda(x^{k_0})}^{k_0}}x^{k_0\rq}\times z$, for some $0 <n_{k_0}+\dots +n_p{-} m_q = \vert z\vert < n_{k_0} {-}j_{\lambda(x^{k_0})}^{k_0}$.
\medskip

Consider the element  $\tilde {x}^{k_0} := x^{k_0}\vert _{\{1,\dots ,n_{k_0}+\dots +n_p{-}m_q\}}\in {\mbox {\bf ST}_{n_{k_0}+\dots +n_p{-}m_q}^l}$, for some $l\geq 1$. We get that
$$f\vert _{\{1,\dots ,r_1+\dots +r_{k_0{-}1}+l\}}\circ (x^1\t \ldots \t {\tilde {x}}^{k_0}) = y^1\cdot \ldots \cdot y^{q-1}.$$
\medskip

So, $q{-}1 \leq k_0 \leq {-}1$ by inductive hypothesis, which implies that if $k_0 < p{-}1$, then $q < p$ and the proof is over for this case.
\medskip

To end the proof, assume that $k_0 = p{-}1$.  As $y^q$ is irreducible, we have that:
$$y^{q\rq} = f\rq \vert _{\{r_1+\dots +r_{p{-}1}{-}p{-}l+2,\dots ,r_1+\dots +r_p{-}p\}} \circ (z\t x^{p\rq}),$$
where $1_{l+r_p-1} < f\rq \vert _{\{r_1+\dots +r_{p{-}1}{-}p{-}l+2,\dots ,r_1+\dots +r_p{-}p\}}$.
\medskip

Moreover, we get that 
$$\prod _{r_1<\dots <r_1+\dots +r_{p-1}} f\rq \vert _{\{1,\dots ,r_1+\dots +r_{p{-}1}{-}l{-}p+1\}} \circ (x^1\t\ldots \t {\tilde {x}}^{p-1}) = y^1\cdot \ldots \cdot y^{q-1}.$$\end{enumerate}
\medskip

So, by recursive hypothesis, $q\leq p$, and we have that $1_{r_1+\dots +r_p{-}p} < f\rq$, which ends the proof.

\end{proo}
\bigskip

\begin{definition} \label{firstreduction} For $x\in {\mbox {\bf Prod}(p)}$ such that $\psi_{\mathbb M}(x) = \sum _{a_i\neq 0}a_ix_i$, define 
$$\psi _{{\mathbb M}1}(x) := \sum _{{x_i\in Prod(p)}\atop {a_i\neq 0}} a_i x_i,$$
where the sum is taken over all the terms $x_i$ appearing in $\psi _{\mathbb M}(x)$ such that $x_i\in {\mbox {\bf Prod}(p)}$.\end{definition}
\medskip

Proposition \ref{proposition4} implies that if the set 
$$\{ \psi _{{\mathbb M}1}(x)\mid x = x^1\cdot \ldots \cdot x^p, x^j\in {\mbox {\bf Br}}\ {\rm for}\ 1\leq j\leq p\}$$
is linearly independent, then $\{ \psi _{{\mathbb M}}(x)\mid x = x^1\cdot \ldots \cdot x^p, x^j\in {\mbox {\bf Br}}\ {\rm for}\ 1\leq j\leq p\}$ is linearly independent, too.
\bigskip

For the next reduction, we need some additional results.

\begin{lemma} \label{lemma9medio} Let $x=\prod _{j_1<\dots < j_{\lambda(x)}}x\rq\in {\mbox {\bf Br}}_n^r$ and let ${\underline l}= \{1\leq l_1 < \dots <l_p < r\}$ be such that $x^{-1}(\{1,\dots , l_p\})\subseteq \{1,\dots ,j_1-1\}$. Given a permutation $f\in {\mbox {Sh}^{\succ}(l_p, r{-}l_p)}$, the element $f\circ (x\vert ^{\{1,\dots ,l\}}\t \dots \t x\vert^{\{l_{p-1}+1,\dots ,l_p\}}\t x\vert^{\{l+1,\dots ,r\}})$ is different to $x$.\end{lemma}
\medskip

\begin{proo} In order to simplify notation, we assume that $p =1$, the proof of the general case is identical. 

For $1\leq l< r-1$, suppose that the set $x^{-1}(\{1,\dots , l\}) = \{k_1<\dots < k_s\}$, for $k_s< j_1-1$. 

Let $h_1$ be the minimal element such that 
$q_1 = x(h_1) = {\mbox {max}\{ x(\{ 1,\dots ,k_s\})\}}$. As is is irreducible, we know that $h_1$ exists and is smaller than $ k_s$, and  there exists at least one $h_2 > k_s$ such that $x(h_2) = q_2 \leq q_1$.
\medskip

Let $f\in {\mbox {Sh}(l, r{-}l)}$. 

\noindent If $f(l) < q_1$, then the minimal element of $(f\circ (x\vert ^{\{1,\dots ,l\}}\t x\vert ^{\{l+1,\dots , r\}}))^{-1}(q_1)$ is $h_1 + m_1$, where $m_1$ is the cardinal of $\{k_i\mid k_i > h_1\}$. 

\noindent So, $f\circ (x\vert ^{\{1,\dots ,l\}}\t x\vert ^{\{l+1,\dots , r\}})\neq x$.
\medskip

For $f(l) \geq q_1$, let $m_2$ be the number of elements $1\leq j < h_2$ such that $x(j) > q_2$, and let $m_3\geq 1$ be the number of elements $1\leq i\leq s$ such that $f(x(k_i))\geq q_1$. 
It is clear that the cardinal of the set $${\mbox {}\{ j\mid 1\leq j < h_2\ {\rm such\ that}\ f\circ (x\vert ^{\{1,\dots ,l\}}\t x\vert ^{\{l+1,\dots , r\}})(j) > f\circ (x\vert ^{\{1,\dots ,l\}}\t x\vert ^{\{l+1,\dots , r\}})(h_2)\}}$$ is greater or equal to $m_2 + m_3$, which implies that $f\circ (x\vert ^{\{1,\dots ,l\}}\t x\vert ^{\{l+1,\dots , r\}})$ is different to $x$.
\end{proo}
\medskip

\begin{theorem} \label{compatwithwBo} Let $x^1,\dots ,x^p$ be a family of elements such that $x^j\in {\mbox {\bf Br}_{n_j}^{r_j}}$, $1\leq j\leq p$, and $x = x^1\cdot \ldots \cdot x^p\in {\mbox {\bf Prod}(p)_{\lambda, {\overline m}}}$. 

\noindent If $\psi_{{\mathbb M}1}(x^1\cdot \ldots \cdot x^p) = \sum _{a_u\neq 0} a_u w_u$, then\begin{enumerate}
\item there does not exist an element $w_u$ such that $w_u <_{wB} x^1\cdot \ldots \cdot x^p$,
\item if $w_u= y^1\cdot \ldots \cdot y^p\in {\mbox {\bf Prod}(p)}$ is a minimal element for $\leq _{wB}$, and $w_u\neq x^1\cdot \ldots \cdot x^p$, then $(\vert y^p\vert ,\dots ,\vert y^1\vert ) < (\vert x^p\vert , \dots ,\vert x^1\vert )$ for the lexicographic order.
\end{enumerate} \end{theorem}
\medskip

\begin{proo} $(1)$ For $p=1$, the result is clear. 
\medskip

Let $\psi_{\mathbb M}(x^k) = \sum _{a_{i_k}\neq 0}a_{i_k}x_{i_k}^k$, for $k=1,\dots , p$. We need to prove that for any collection of elements 
$x_{i_1}^1,\dots ,x_{i_p}^p$ and any $f\in {\mbox {Sh}^{\bullet}(r_1,\dots ,r_p)}$, we have that $f\circ (x_{i_1}^1\t \dots \t x_{i_p}^p) \not < _{wB} x^1\cdot \dots \cdot x^p$.

The permutation $f\rq \in {\mbox {Sh}(r_1{-}1,\dots ,r_p{-}1)}$ may be written in a unique way as
$$f\rq = f_2\circ (f_1\t 1_{r_p{-}1})\in {\mbox {Sh}(r_1{-}1,\dots ,r_p{-}1)},$$
where $f_1\in {\mbox {Sh}(r_1{-}1,\dots ,r_{p{-}1}{-}1)}$ and 
 $f_2\in {\mbox {Sh}(r_1+\dots +r_{p{-}1} {-}p+1,r_p{-}1)}$. 

Applying recursive hypothesis, we know that 

$$\prod _{r_1<\dots < r_1+\dots +r_{p{-}1}}f_1\circ (x_{i_1}^1\t \dots \t x_{i_{p{-}1}}^{p{-}1}) \not < _{wB} x^1\cdot \dots \cdot x^{p{-}1}.$$ 

So, it suffices to prove the result assuming that $x_1\in {\mbox {\bf Prod}(p{-}1)}$ and $x_2\in {\mbox {\bf Br}}$.
\medskip

Applying Remark \ref{elementsPL}, the element $w_u= f\circ (x_{i_1}^1\t x_{i_2}^2)$ satisfies that 
$x_{i_1}^1 = \prod _{j_1^1<\dots < j_{\lambda (x^1)}^1}(w^1\t z^1)$ and $x_{i_2}^2 = \prod _{j_1^2<\dots < j_{\lambda (x^2)}^2}(y^2 \t w^2)$, 

\noindent where $\vert z^1\vert < n_1{-}j_{\lambda (x^1)}$, $\vert y^2\vert < j_1^2 {-}1$ and $f\in {\mbox {Sh}^{\bullet}(r_1,r_2)}$. 
\medskip

We want to prove that in any case, whenever $f\circ (x_{i_1}^1\t x_{i_2}^2)= {\overline {x}^1}\cdot {\overline {x}^2}$, for ${\overline {x}^1}$ in ${\mbox {\bf Prod}(p{-}1)}$ and ${\overline {x}^2}\in {\mbox {\bf Br}}$ , we have $g\circ ({\overline {x}^1}\cdot {\overline {x}^2}) \neq x^1\cdot x^2$.
\medskip

For elements $f =\prod _{r_1< r_1+r_2} f\rq$,  $w^1 \in {\mbox {\bf ST}_{m_1}^{s_1}}$, $z^1\in {\mbox {\bf ST}_{m_1}^{r_1{-}s_1{-}1}}$ and $y^2 \in {\mbox {\bf ST}_{m_2}^{s_2}}$, we have that. 
\begin{enumerate}
\item If $f\rq (s_1) > f\rq (r_1 +s_2)$, then $f\circ (x_{i_1}^1\t x_{i_2}^2)$ belongs to ${\mbox {\bf Prod}(p{-}1)}$. 
\item If $f\rq (r_1) < f\rq (s_1)$ and $f\rq (r_1{-}1) < f\rq (r_1+s_2)$, then:\begin{enumerate}
\item ${\overline {x}^1} = \prod _{j_1^1< \dots < j_{\lambda (x^1)}^1}( f_1\rq\circ (w^1 \t z^1 \t y_1^2)\times y_2^2)$, with $y^2 = y_1^2\times y_2^2$  and $f_1\rq \in {\mbox {Sh}(r_1{-}1,l_1)}$,
\item ${\overline {x}^2} = \prod _{j_1^2< \dots < j_{\lambda (x^2)}^2}w^2 $.\end{enumerate}
\item If $f\rq (r_1)< f\rq (s_1)$ and $f\rq (r_1+s_2) < f\rq (r_1-1)$, then $f\circ (x_{i_1}^1\t x_{i_2}^2)$ belongs to ${\mbox {\bf Prod}(p{-}1)}$. 
\item If $f\rq (s_1) < f\rq (r_1)$ and $f\rq (r_1+s_2) < f\rq (r_1-1)$, then \begin{enumerate}
\item ${\overline {x}^1} = \prod _{j_1^1< \dots < j_{\lambda (x^1)}^1}(w^1 \t z_1^1)$, with $z^1 = z_1^1\times z_2^1$,  
\item ${\overline {x}^2} = \prod _{j_1^2< \dots < j_{\lambda (x^2)}^2}f_2\rq \circ (z_2^1\t y^2\t w^2 )$, with $f_2\rq \in {\mbox {Sh}(l_2,r_2{-}1)}$.\end{enumerate}
\item If $f\rq (r_1{-}1) < f\rq (r_1+s_1)$, then \begin{enumerate}
\item ${\overline {x}^1} = \prod _{j_1^1< \dots < j_{\lambda (x^1)}^1}(w^1 \t f_3\rq \circ(z^1\t y^2))$, with $f_3\rq \in {\mbox {Sh}(r_1{-}s_1{-}1,s_2)}$,  
\item ${\overline {x}^2} = \prod _{j_1^2< \dots < j_{\lambda (x^2)}^2}w^2 $.\end{enumerate}
\end{enumerate}

The unique cases where $f\circ (x_{i_1}^1\t x_{i_2}^2)$ is the product $\cdot$ of $p$ elements are $(2),\ (4)$ and $(5)$.
\medskip

Suppose now that we have ${\overline {x}^k} = \prod _{j_1^k< \dots < j_{\lambda (x^k)}^k} {\overline {x}^k} \rq \in {\mbox {\bf ST}_{q_k}^{v_k}}$, for $k = 1,2$, and $g\in {\mbox {Sh}^{\bullet}(v_1,v_2)}$.

If $\vert {\overline {x}^1} \vert < n_1$, then ${\overline {x}^2}$ must be reducible in order to get $g\circ ({\overline {x}^1}\t {\overline {x}^2}) = x^1\cdot x^2$ but ${\overline {x}^2}$ is irreducible, so in the case $(4)$ there is no solution.
\medskip

In case $(2)$, we have that if $g\circ ({\overline {x}^1}\t {\overline {x}^2}) = x^1\cdot x^2$, then $n_2 \leq n_2 {-} m_2 + \vert y_2^2\vert $, but $\vert y_2^2\vert < m_2$, and there does not exist a solution.
\medskip

In $(5)$, we need that $f_3\rq\circ  (z^1\t y^2)) = t^1\t t^2$, with $\vert t^2\vert = \vert y^2\vert = m_2$. But if $f_3\rq $ is not the identity, then $f_3\rq\circ  (z^1\t y^2)$ cannot be decomposed as $ t^1\t t^2$, with $\vert t^2\vert = m_2$. 

If $f_3\rq $ is the identity, then we need that $x_2 = g_1\circ (y_2\times {\overline {x}_2})$. 

As $x^2 = y\backslash w$ is in ${\mbox {\bf Br}}$, there exists ${\underline l} = (1\leq l_1< \dots < l_p)$ such that $y^2 = y\vert ^{\{1,\dots ,l_1\}}\t \dots \t  y\vert ^{\{l_{p-1},\dots ,l_p\}}$. But, from Lemma \ref{lemma9medio}, we get that there does not exist $g_1$ such that $x^2 = g_1\circ (y^2\times {\overline {x}^2})$, which ends the proof of point $(1)$.
\bigskip

\noindent $(2)$ Applying the same argument that in point $(1)$, it suffices to prove the assertion for $x^1\in {\mbox {\bf Prod}(p{-}1)}$ and $x^2\in {\mbox {\bf Br}}$.

Suppose that $x^1$ and $x^2$, are such that $\psi_{\mathbb M}(x^1) = \sum _{a_i\neq 0} a_i x_i^1$ and $\psi _{\mathbb M}(x2) = \sum _{b_j\neq 0} a_j x_j^2$. 

We have that $w_u = f\circ (x_l^1\t x_k^2) = w_1\cdot w_2$, for some pair $l,k$ and $f\in {\mbox {Sh}^{\bullet}(r_1,r_2)}$.

Let \begin{enumerate}
\item $x_l^1 = \prod _{j_1<\dots <j_{\lambda (x_1)}} (x_l\rq \t z_1)$, for $0 < \vert z_1\vert < n_1-j_{\lambda (x_1)}$, 
\item $x_k^2 =  \prod _{h_1<\dots < h_{\lambda (x_2)}} (y_2\t x_k\rq )$, for $0 \leq \vert y_2\vert < h_1-1$. \end{enumerate}

As $\vert x_2^k\vert = \vert x_2\vert$, if $\vert x_2^k\vert \leq  \vert w_2\vert $, then 
$$w_2 = f\vert _{\{k+1, \dots ,r_1+r_2{-}2\}}\circ (z_{12} \t x_k^{2\rq} ),$$
where $z_1= z_{11}\t z_{12}$, $k+ 1 \leq r{-}1$ and $f\rq \in {\mbox {Sh}(r_1{-}k{-}1, r_2-1)}$.

As $w_2$ is irreducible, we get that $f\rq \neq 1_{r_1+r_2+{-}k{-}2}$ is not the identity, and therefore $x_l^11\cdot x_k^2 <_{wB} w_1\cdot w_2$. So, $w_u$ is not minimal.
\end{proo}
\bigskip

Applying Theorem \ref{compatwithwBo}, let 
$$\psi_{Prod}(x_1\cdot \ldots \cdot x_p) := \sum _{{w_u\ {\rm minimal\ for}\ \leq_{wB}}\atop {a_u\neq 0}} a_u w_u,$$
where $\psi_{{\mathbb M}1}(x_1\cdot \ldots \cdot x_p) = \sum _{a_u\neq 0}a_u w_u$. We have that:\begin{enumerate}
\item ${\psi}_{Prod}(x_1\cdot \ldots \cdot x_p) = x_1\cdot \ldots \cdot x_p + \dots $,
\item if $w_u = y_1\cdot \ldots \cdot y_p\in {\mbox {Prod}(p)}$ is minimal and different from $x_1\cdot \ldots \cdot x_p$, then 
$(\vert y_p\vert ,\dots ,\vert y_1\vert ) < (\vert x_p\vert ,\dots ,\vert x_1\vert )$ for the lexicographic order,
\item if the set $\{ \psi_{Prod}(x)\mid x\in {\mbox {\bf Prod}}\}$ is linearly independent, then the $\{ \psi _{\mathbb M}(x)\mid x\in {\mbox {\bf Prod}}\}$ is linearly independent, too.\end{enumerate}

\bigskip

Finally, applying Proposition \ref{proposition4} to $x = x^1\cdot \ldots \cdot x^p\in {\mbox {\bf Prod}(p)}$, we get that:
$$\psi_{Prod}(x) = \sum _{a_i\neq 0} a_i w_i,$$ 
with $w_i = y_i^1\cdot \ldots \cdot y_i^p\in {\mbox {\bf Prod}(p)}$ such that $(\vert y_i^p\vert, \dots ,\vert y_i^1\vert ) < (\vert x^p\vert , \dots , \vert x^p\vert )$ for the lexicographic order. 
So, if the set $\{ x^1\cdot \ldots \cdot x^p\mid x^1,\dots ,x^p\in {\mbox {\bf Br}}\}$ is linearly independent for any integer $p\geq 2$, we get that
$$\{ \psi ^0 (x^1\cdot \ldots \cdot x^p)\mid x^1,\dots ,x^p\in {\mbox {\bf Br}},\ p\geq 2\}$$
is linearly independent, too.
\bigskip

\bigskip

\centerline{\bf Proof of Theorem \ref{principthe}}
\medskip

In the previous section, we have shown that we may restrict ourselves to prove that $\psi ^0: \K[{\mbox {\bf Irr}}]\longrightarrow {\mbox {Prim}(\st_{qT})}$ is an isomorphism.

We have that ${\mbox {\bf Irr}}$ is the disjoint union of $\Br$ and $\Prod$. Moreover, we have that $E({\mbox {\bf Irr}}) = {\mbox {Prim} ({\st})}$ by Remark \ref{idempotentE}. Let us describe the framework of our proof.
\medskip

As $E(x) - x$ belongs to $\K[{\mbox {\bf Red}}]$, for all $x\in {\mathcal B}$, it is immediate that $\{\psi _{\mathbb M}(x)\mid \ x\in {\mathcal B}\}$ is linearly independent in ${\mbox {Prim}(\st)}$.
\medskip

 Corollary \ref{corollarybrace} shows that $\psi_{\mathbb M}({\mbox {\bf Br}})\subseteq \K [{\mbox {\bf Br}}]$. We shall prove first that $\{\psi _{Br}(x)\mid \ x\in {\mbox {\bf Br}}\}$ spans $\K[{\mbox {\bf Br}}]$, which implies that both sets $\{\psi _{\mathbb M}(x)\mid \ x\in {\mbox {\bf Br}}\}$ and $\{\psi ^0(x)\mid \ x\in {\mbox {\bf Br}}\}$ are linearly independent 
 in $\K [{\mbox {\bf Br}}]$.
 \medskip

Using the result above, it is immediate to see that the set 
$$\psi_{\mathbb M}({\mbox {\bf Br}})_n^{\bullet}: = \{ \psi _{\mathbb M}(x_1)\cdot \ldots \cdot \psi_{\mathbb M}(x_p)\mid x_j\in {\mbox {\bf Br}_{n_j}},\ 1\leq p\leq n\ {\rm and}\ \sum _{j=1}^p n_j =n\},$$
is linearly independent. To end the proof, we show that the spaces spanned by the sets $\psi_{\mathbb M}({\mbox {\bf Br}})_n^{\bullet}$ and 
$\psi_{\mathbb M}({\mbox {\bf Irr}_n})$ are equal, for all $n\geq 1$. We proceed by induction on $n$.
\bigskip

For $n=1$ and $n=2$, the result is evident.

For $n\geq 3$, from Corollary \ref{corollarybrace}, we know that if $\{\psi_{Br}(x)\mid x\in \Br_m \}$ and $\{\psi ^0(x)\mid x\in \Prod _m \}$ are linearly independent for all $m <n$, hence the set $\{\psi ^0(x)\mid x\in {\mbox {\bf Irr}_m}\}$ is linearly independent in $\K[{\mbox {\bf ST}_m}]$, for all $m < n$.

So, the image of the subspace ${\displaystyle \bigoplus _{m=1}^{n-1}\K[{\mbox {\bf Irr}_m}]}$ under $\psi^0$ spans the subspace ${\displaystyle \bigoplus _{m=0}^{n-1}{\mbox {Prim}(\st)_m}}$.
\medskip

\subsection{Step $\Br$} We proceed as in Section \ref{dendbialg}. 

By induction, we assume that $\{ \psi ^0(x) \mid  x\in {\mbox {\bf Irr}_m},\ {\rm for}\ m<n\}$ spans ${\mbox {Prim} (\st)_{<n}}:={\displaystyle \bigoplus _{j=1}^{n-1}{\mbox {Prim} (\st)_j}}$.

Applying Proposition \ref{propositionbrace} we get that 
$$\bigoplus _{j=1}^{n-1}\K[{\mbox {\bf ST}_j}]\subseteq \bigoplus _{j=1}^{n-1}\omega ^{\succ}\bigl ({\mbox {Prim} (\st)_{<n}}^{\ot j}\bigr ).$$

So, the set ${\displaystyle \bigcup_{x\in \B_{\leq n}}\{ x\backslash (\omega ^{\succ} \bigl ({\mbox {Prim} (\st)_{<n}}^{\ot j}\bigr ))\}}$ spans ${\displaystyle \bigcup _{j=0}^n\K[\Br _j]}$, which implies that $\{\psi_{Br}(x)\mid x\in \Br\ {\rm and}\ \vert x\vert \leq n\}$ is linearly independent in $\K[{\mbox {\bf ST}}]$. Therefore, the set $\{\psi ^0(x)\mid  x\in \Br\ {\rm and}\ \vert x\vert \leq n\}$ is linearly independent in ${\mbox {Prim}(\st)}$.
\bigskip, 

\subsection{Step $\Prod$} 
\medskip

We need to prove that for any $p\geq 2$, the set 
$$\{ x^1\cdot \ldots \cdot x^p\mid x^1,\dots ,x^p \in {\mbox {\bf Br}}\}$$ is linearly independent, which is obviously true.

\bigskip

\bigskip

\subsection*{Final comment} In order to simplify notation, our definitions of $q$-dendriform algebra, brace algebra and ${\mbox {GV}_q}$ algebra were given in the non-graded case. A graded version of these notions is obtained just applying the Koszul sign convention, and our results still hold. In particular, F. Chapoton\rq s operad of $K$-algebras in \cite{Chapoton}, which is  described by permutohedra, coincides with a differential graded version of $0$-tridendriform algebras, such that the degree of the operations $\succ $ and $\prec $ is $1$ while the degree of the associative product $\cdot $ is $0$. As $0$-tridendriform is non symmetric, the unique relation which is modified by Koszul\rq s sign in Definition \ref{def:TriDend} is:
$$x\cdot (y\succ z) = (-1) ^{\vert y\vert } (x\prec y)\succ z.$$

In this case, the coboundary map of the permutohedra is described as the unique differental map $\partial$  of degree $-1$ such that:\begin{enumerate}
\item $\partial (x\succ y) = \partial (x)\succ y + (-1)^{\vert x\vert } x\succ \partial (y) - (-1)^{\vert x\vert +1} x\cdot y$,
\item $\partial (x\cdot y) = \partial (x)\cdot  y + (-1)^{\vert x\vert +1} x\cdot \partial (y)$,
\item $\partial (x\prec y) = \partial (x)\prec y + (-1)^{\vert x\vert } x\succ \partial (y) - (-1)^{\vert x\vert } x\cdot y$,\end{enumerate}
for homogeneous elements $x,y,z\in \st$.

In this case, our result implies that the ${\mbox {GV}_0}$ algebra structure on ${\mbox {Prim}(\st)}$ equipped with the coboundary map of the permutohedra, is a free cacti algebra (see for instance \cite{Kaufmann}) on the base ${\mathcal B}$.

\end{document}